\documentclass[red,11pt,a4paper]{article}

\usepackage{cite}

\usepackage{amsmath}
\usepackage{amscd}
\usepackage{amssymb}
\usepackage[pdftex]{color,graphicx,hyperref}
\usepackage{latexsym}
\usepackage{color}
\usepackage[normalem]{ulem}
\bf
\usepackage{ika}

\title{Two-velocity hydrodynamics in fluid mechanics: Part II \\
        Existence of global $\kappa$--entropy solutions \\
         to compressible Navier-Stokes systems 
        \\ with degenerate viscosities}
\author{Didier Bresch$^1$,  Beno\^{\i}t Desjardins$^{2}$, Ewelina Zatorska$^{3,4,5}$}

\begin{document}
\maketitle
\normalsize
\begin{center}
{\small 1.  Universit\'e de Savoie, Laboratoire de Math\'ematiques\\
73376 Le Bourget du Lac, France\\
      {\ }\\
2.  Fondation Math\'ematique Jacques Hadamard and  CMLA ENS Cachan\\ 
	94235 Cachan cedex\\ {\ }\\
3.  Centre de Math\'{e}matiques Appliqu\'{e}es,\\
      \'{E}cole Polytechnique, 91128 Palaiseau Cedex, France\\
      {\ }\\
4. Institute of Mathematics\\
Polish Academy of Sciences, ul \'Sniadeckich 8, 00-656 Warszawa, Poland\\
      {\ }\\
5. Institute of Applied Mathematics and Mechanics\\
University of Warsaw, ul. Banacha 2, 02-097 Warszawa, Poland\\}
%{\sc E-mail: e.zatorska@mimuw.edu.pl}}
\end{center}

\noindent{\bf Abstract:} This paper addresses the issue of global existence of so-called $\kappa$--entropy solutions to the Navier--Stokes equations for viscous compressible and  barotropic fluids with degenerate viscosities. We consider the three dimensional space domain with periodic boundary conditions.  Our solutions satisfy the weak formulation of the mass and momentum conservation equations and also a generalization of the BD--entropy identity called: {\it $\kappa$--entropy}.  This new entropy involves  a mixture parameter $\kappa \in (0,1)$  between the two velocities $\vu$ and $\vu+2\Grad\varphi(\vr)$ (the latter was introduced by the first two authors in [C.~R. Acad. Sci. Paris 2004]), where $\vu$ is the velocity field and $\varphi$ is a function of the density $\vr$ defined by $\varphi'(s)=\mu'(s)/s$. 
%The assumption $\lambda(\vr)=2(\mu'(\vr)\vr -\mu(\vr))$ is also required as  in previous works by the two first authors to establish the $\kappa$-entropy identity.
  As a byproduct of the existence proof,  we show that two-velocity hydrodynamics  (in the spirit of {\sc S.C. Shugrin} 1994) is a possible formulation of a model of barotropic compressible flow with degenerate viscosities. It is used in construction of approximate solutions based on a conservative  augmented approximate scheme.
  %%%% I don't know what is the information following from this part below%%%%%%
%  of the compressible Navier-Stokes equations extending  approach 
%  developed for some inviscid compressible systems with capillarity by other authors.

\medskip

\noindent {\bf Keywords.} Hypocoercivity, Compressible Navier-Stokes, Augmented system,
Two-velocity hydrodynamics, $\kappa$-entropy.

\section{Introduction}
In 2006--2007, the first two authors introduced the concept of global weak solutions to the Navier--Stokes equations for compressible barotropic fluids degenerate viscosities satisfying the energy inequality and  an extra mathematical entropy called BD entropy, \cite{BrDe4}, \cite{BrDe2}, \cite{BrDe3}, 
  The BD entropy identity was derived in \cite{BrDe1} for $\lambda(\vr)= 2(\mu'(\vr)\vr - \mu(\vr))$ as a generalization of  the work of the first two authors and C.K. {\sc Lin} \cite{BrDeLi}, where  $\mu(\vr)=\vr$ and $\lambda(\vr)=0$. It involves an energy related to the velocity $\vu+ 2\Grad\varphi(\vr)$ where $\varphi$ is a function of the density $\vr$ defined by 
  $\varphi'(s)=\mu'(s)/s$. It is worth to note that such quantity appears also in study of inviscid systems,  see \cite{BeDaGrSaSm} where $\varphi(\vr)=\log\vr$ and in the works by E. {\sc Nelson} (see for example \cite{Ne}) related to kinematics of Markovian motion with interesting discussion using two velocities formulation. There, quantity $\mu \nabla\vr/\vr$ stands for velocity required for the particle to counteract osmotic effects (osmotic velocity) and the current velocity.

Including a drag force in the momentum equation or an additional singular pressure terms, stability of global weak solutions to the barotropic compressible Navier--Stokes equations with density dependent viscosities satisfying the extra BD entropy was proven in \cite{BrDe3}, \cite{BrDe4} and with D. {\sc G\'erard-Varet} in \cite{BrDeGe}, see also 
 %with no need of multiplying the momentum equation by $\vr$ as in 
 \cite{BrDeLi}. 
Stability of global weak solutions satisfying the extra BD entropy without these extra terms was obtained by {\sc A. Mellet} and {\sc A. Vasseur} in \cite{MeVa}. They showed additional estimate on the velocity field but the construction of approximate solutions satisfying the energy, the BD entropy and the Mellet--Vasseur estimate together was so far an open problem. Very recently {\sc A. Vasseur}  and {\sc C.~Yu} have proved in \cite{VaYu} that given a global  weak solution to the compressible Navier--Stokes equations with turbulent drag terms and appropriate capillary term satisfying the energy and the BD entropy estimates, it is possible to pass to the limit and get a weak solution to the compressible Navier-Stokes system. They constructed smooth multipliers allowing to get the Mellet--Vasseur estimate uniformly with respect to the drag coefficient and to pass to the limit in the capillary terms.
 The uniform control allows to suppress the drag terms letting the drag coefficients tend to zero.
 This result, when coupled with our present paper gives the first complete existence result to the compressible Navier--Stokes equations with general degenerate viscosities with no extra terms.
 Here we construct the weak solution to the system similar to the considered in \cite{VaYu} modulo capillarity terms. It was however shown in \cite{BrGiZa}  that this construction is compatible with the  quantum capillary term  of the form $\kappa \vr \nabla (\Delta \sqrt\vr / \sqrt\vr)$ appearing in the ghost system \cite{LeSuTr}, see also  \cite{GiVi} for the study on quantum viscous Navier-Stokes system and \cite{BrCoNoVi} for a full range of compatible capillary terms.

    Concerning  construction of approximate solutions with singular pressure, drag terms or capillarity terms,  the authors gave some hints in \cite{BrDe2} for the general setting $\lambda(\vr)= 2(\mu'(\vr)\vr - \mu(\vr))$.  This has been fully developped in \cite{EZ2}, \cite{MuPoZa1}, \cite{MuPoZa2} in the case of a linear viscosity $\mu(\vr)=\vr$. However, even for this case, the procedure seems to be very complex. The approximation involves: i) regularization of the mass equation, ii) an extra regularizing terms in the momentum equation inspired by those given in \cite{BrDe2}, iii) high-order generalized capillarity term,  and iv) regularization of the velocity $\vu + 2\Grad\varphi(\vr)$.  
   
In the present paper, we propose to come back to the compressible Navier-Stokes equations with
degenerate viscosities in the general setting with simpler construction of approximate solutions.
   We propose a very natural concept of global $\kappa$--entropy solutions based on a generalized BD entropy. The global weak solution of the compressible Navier-Stokes equations which satisfies the BD entropy is also a $\kappa$--entropy solution for all $0<  \kappa < 1$, however, the starting point for the proofs of existence are rather different. 
  Through $\kappa$-entropy  we introduce a two-velocity hydrodynamical formulation of the compressible Navier--Stokes equations with degenerate viscosities.
  Our construction of the approximate solutions is  based on an augmented approximate scheme using this two-velocity structure.  Similar approximate scheme has been introduced  for the zero Mach number system in  the first part of the present series \cite{BrGiZa}.  For reader's convenience, since lack of divergence-free condition gives rise to several new terms and since the $\kappa$-entropy and the limit process changes,  we repeat the whole construction from \cite{BrGiZa}.
  
    By using augmented systems, we extend an approach developed for inviscid compressible systems with dispersion (see for instance {\sc S. Benzoni, R. Danchin, S. Descombes}  \cite{BeDaDe} and {\sc F. B\'ethuel, R. Danchin, P. Gravejat, J.--C.~Saut, D. Smets}  \cite{BeDaGrSaSm}). 
  Lately, similar structure was also used by {\sc P.~Noble} and {\sc J.--P.~Vila} \cite{NoVi} to study the stability of various  approximations of the one-dimension Euler--Korteweg equations (dispersive system). They introduced an additional unknown (the gradient of a function of the density) in order to rewrite the system into a hyperbolic system perturbed by a second order skew symmetric term. 
Derivation of relevant  numerical scheme using an additional unknown for  Navier-Stokes-Korteweg system (dispersive system) is purpose of a forthcoming paper \cite{BrNoViVi}.

In our construction, the only regularization of the continuity equation  is via regularization of the new velocity $\vw=\vu+2\kappa\Grad\varphi(\vr)$ and no  further viscous approximation is needed.  Parabolicity of the equation for the density written in this form has been recently observed and studied by {\sc B. Haspot} in \cite{Ha14} with the change of variables $\vv=\vu + \nabla\varphi (\vr)$. This property was then used to obtain result concerning global controllability for the shallow-water system with two control forces, see \cite{DrHa}.

     The change of unknowns from the work of  {\sc B. Haspot} corresponds to $\kappa=1/2$ from our definition of $\kappa$-entropy solution.  The 1/2--entropy solution has been obtained recently by {\sc M. Gisclon} and {\sc I.~Violet} \cite{GiVi}, for $\mu(\vr)=\vr$ and $\lambda(\vr)=0$, starting from the quantum compressible Navier-Stokes equations studied by {\sc A. J\"uengel} in   \cite{Ju} (extended in \cite{Do}  and \cite{Ji}) with an additional singular pressure and letting the scaled Planck constant vanish. Their proof strongly relies on the Bohm potential identity and therefore works only for $\mu(\vr)=\vr$ in the multi-dimensional space case.  See also result by {\sc B. Haspot} concerning Korteweg systems in \cite{Ha2} and references cited therein. 
     
Interestingly enough, it turns out that our two-velocity formulation is linked to \cite{Sh} and \cite{GaSh}: we introduce a
 mixture parameter $\kappa \in (0,1)$ which combines two velocity vector fields $\vu$ and  $\vu + 2\Grad\varphi(\vr)$ sharing the same reference density. 
   The augmented system used to construct approximate solution is also similar to the one from \cite{Sh} where the main velocity $\vu+ 2\kappa\Grad\varphi(\vr)$ and a drift  $2\Grad\varphi(\vr)$ correspond to the same reference density $\vr$. The reader is also  referred to the paper by E. {\sc Feireisl} and {\sc A. Vasseur} \cite{FeVa} where they study a  compressible system with two-velocities proposed by {\sc H.  Brenner}. 
        
     The generalized BD entropy ($\kappa$-entropy) reflects some non-linear hypocoercivity property of the nonlinear compressible Navier-Stokes equations under the aforementioned relation between $\lambda$ and $\mu$. For an introduction to hypercoercivity, the  interested reader is referred to the paper by {\sc K. Beauchard} and {\sc E. Zuazua} \cite{BeZu} and  the book  by C.~{\sc Villani} \cite{Vi}  (and references cited therein) which describe its link to global existence around equilibrium and large-time behavior (see  \cite{Da} by {\sc R. Danchin} for an application to fluid mechanics). Readers interested in entropy for nonlinear partial differential equations are also referred to \cite{BaCh}.
    Our generalized BD entropy ($\kappa$-entropies) may be seen as a nonlinear version of the identity that was proven by {\sc A. Matsumura--T. Nishida} on the linearized compressible system around the equilibrium $(\vr_{\rm eq},\vu_{\rm eq})=(1,0)$. For the reader's convenience, we revisit hypocoercivity on  linearized  compressible Navier--Stokes systems  in the last two sections of the paper (barotropic and heat-conducting case). 
%    This shows that $\kappa$-entropy estimate for such linearized system with a coefficient $\kappa$  may be chosen arbitrarily small.
%   The definition of the global $\kappa$--entropy solution is strongly associated with two velocity hydrodynamics in the compressible Navier--Stokes equations with degenerate viscosities. 
%   The reader interested  on the subject is referred for instance to \cite{Sh} for some discussions
%   around two velocity Hydrodynamics.

%              In some sense, it seems that our result is the first showing clearly the existence of two-velocity hydrodynamics  in the  barotropic compressible Navier--Stokes equations with degenerate viscosities. 
%   The augmented system is closed to systems written in \cite{Sh} (with a joint  density  for the two species) with the main velocity given by $\vw= \vu+2\kappa\Grad \vp(\vr)$ and the drift by $\vv= 2\Grad \varphi(\vr)$.

The purpose of the last section is also to present a thermodynamically consistent two-velocity model with heat conductivity in the spirit of the work of {\sc S.M. Shugrin},  \cite{Sh}. Each of the velocity vector fields $\vu$ and $2\Grad\varphi(\vr)$ is associated with different density $(1-\kappa)\vr$ and $\kappa\vr$, respectively. The main objective is not to prove global existence of solutions rigorously but  to show that the two-velocity hydrodynamics is consistent with the study  performed in the first part of the present series  \cite{BrGiZa} for the zero Mach number system. Our formulation uses the generalized $\kappa$-temperature which is not {\it a priori} the usual temperature. In this sense, our formulation and the usual heat-conducting compressible Navier-Stokes equations are not equivalent. However,  the formal low-Mach number limit for such two-velocity system with generalized  temperature gives exactly the augmented system from \cite{BrGiZa} used to obtain the existence result.  This is a kind of consistency between the approaches and the systems.

\section{The barotropic Navier-Stokes system}

Two compressible fluid models with degenerate viscosity and pressure depending only on the density (barotropic flows) will be considered.

\medskip

\noindent {\bf1)  Compressible Navier-Stokes with singular pressure.}  First  the compressible Navier--Stokes equations for compressible and barotropic fluids write as follows:
        \begin{equation}\label{main}
     \left\{ 
      \begin{array}{l}
            \vspace{0.2cm}
       \pt\vr+\Div\lr{\vr\vu}= 0,\\
            \vspace{0.2cm}
        \pt\lr{\vr\vu} + \Div ({\vr\vu\otimes\vu})  - \Div\lr{2 \mu(\vr) D(\vu)} 
           - \Grad (\lambda(\vr) \Div \, \vu) +   \Grad p(\vr) = 0,\\%& \mbox{w}& (0,T)\times\Omega,\\
           \end{array}\right. \quad\mbox{in}\ (0,T)\times\Omega
    \end{equation}
where $D(\vu)=\frac{1}{2}\left(\Grad \vu+\Grad^t \vu\right)$, $(\Div ({\vr\vu\otimes\vu})_j = \partial_i(\vr\vu_i u_j)$ and  $\Omega$ is a periodic box $\Omega=\mathbb{T}^3$.
   The  pressure $p$  is  singular close to zero density
	\begin{equation}\label{coldp}
	p'(\vr)=\left\{
	\begin{array}{cl}
	c_1\vr^{-\gamma^- -1}&\mbox{for} \ \vr\leq \vr^*,\\
	c_2\vr^{\gamma^+-1}&\mbox{for} \ \vr>\vr^*
	\end{array}\right.
	\end{equation}
with $c_i>0$ , $\gamma^+>1$ and $\gamma^- >0$. A more precise estimation of  $ \gamma^-$ 
will be given below. 

\medskip

The viscosity coefficients $\mu(\vr)$, $\lambda(\vr)$ satisfy the  {\sc Bresch-Desjardins} relation introduced in \cite{BrDe1}
	\eq
		{\lambda(\vr)=2(\vr\mu'(\vr)-\mu(\vr)).\label{BD}}
    This system is completed with initial data
\begin{equation}\label{ini}
\vr\vert_{t=0} = \vr_0, \qquad (\vr\vu)\vert_{t=0}= m_0.
\end{equation}

\medskip

The pressure \eqref{coldp} was introduced in \cite{BrDe3}. Note that the singular part is active only for density close to vacuum and that the usual power law equation of state is recovered far from vacuum.
The physical relevance of the compressible Navier-Stokes equations is very questionable in regions where density is close to zero:  the medium is not only unlikely to be in a liquid or gas state (elasticity and plasticity has to be considered for such solid materials, for which by the low densities may lead to negative pressures), but also the rarefied regime of vanishing densities violates the assumptions on the mean free path of particles suitable for fluid models. 
   The negativity of the cold pressure implied by the above assumptions for low densities may also be interpreted as some artificial way to get close to a solid state in tension.  This is exactly the idea
in the work \cite{KiLaNi} where intermolecular forces, namely long-range attractive van der Waals and short-range Born repulsive intermolecular forces, are considered.

From the mathematical point of view, the singular pressure \eqref{coldp} has some stabilizing properties. They  were used in \cite{MuPoZa1} to investigate the model of compressible mixture. Steady compressible Navier-Stokes system with pressure singular at vacuum was also investigated by {\sc M. \L asica} {\rm  \cite{Lasica}}. 
For studies on doubly singular pressure  in the context of mixtures we refer to paper by {\sc E. Feireisl}, {\sc Y. Lu} and {\sc J. M\'alek} {\rm \cite{FeLuMa}}. See also the work by {\sc G. Kitavsev, P. Lauren\c cot, B. Niethammer} {\rm  \cite{KiLaNi}} where they study lubrication equations in the presence of strong slippage.

\bigskip

\noindent {\bf 2)  Compressible Navier-Stokes with a turbulent drag term.} 
We discuss,  in Section \ref{drag}, how to handle the compressible Navier--Stokes equations with turbulent drag term ($r_1\vr|\vu|\vu$, $r_1>0$) and standard gamma-type pressure law, namely
      \begin{equation}\label{main2}
     \left\{ 
      \begin{array}{l}
            \vspace{0.2cm}
       \pt\vr+\Div\lr{\vr\vu}= 0,\\
            \vspace{0.2cm}
        \pt\lr{\vr\vu} + \Div ({\vr\vu\otimes\vu})  - \Div\lr{2 \mu(\vr) D(\vu)} 
           - \Grad (\lambda(\vr) \Div \, \vu) +  r_1 \vr |\vu|\vu +  \Grad p(\vr) = 0,\\
                      \end{array}\right. 
    \end{equation}
 with the usual pressure law $p(\vr)= a\vr^\gamma$ with $\gamma>1$.  
  It is important to cover such situation, see the footnote in the introduction.
  The assumptions on the viscosity in this section will be similar to the one introduced 
 in \cite{MeVa} with the extra control  $\mu(\vr) \le c \vr^{2/3+1/3\nu}$ for $\nu \in (1/\gamma, 1)$ when $\vr \ge 1.$
 This additional assumption is introduced because the energy and BD entropy are mixed in the
 $\kappa$-entropy (less information is available compared to \cite{BrDeGe}).
 In this case, we will only focus on the extra terms that have to be included in the augmented
 system and show how the $\kappa$-entropy changes.  The existence result and
 asymptotic limit with respect to smoothing parameters $\alpha$, $n$ and $\delta$ will not be affected
 when $\varepsilon$ is fixed.  The asymptotic limit when $\varepsilon$  tends to $0$ is then the 
 same as in \cite{BrDeGe}.

%%%%%I think all of this was explained in the introduction%%%%%%%%
%\begin{rmk} The present work may be extended to dispersive compressible 
%Navier-Stokes system (Korteweg type system) when the additional dispersive term is compatible with the velocity  $\vu+2\kappa\Grad\varphi(\vr)$.  
%For instance, when $\varphi(\vr) = \log\vr$,  the quantum dispersive term $- \vr \Grad \Bigl({\Delta\sqrt \vr}/{\sqrt\vr}\Bigr)$ studied in {\rm \cite{Ju}}: see recent works by  {\sc M.~Gisclon}--{\sc I. Violet} or A. {\sc J\"uengel}. Such extension has been discussed
%in the first part {\rm \cite{BrGiZa}} discussing the ghost compressible system.
% See also some interesting result by {\sc B. Haspot} concerning Korteweg systems in {\rm  \cite{Ha2}}.
%  This will be considered in {\rm {\cite{BrCoNoVi}}} and  {\rm \cite{BrNoViVi}}  for numerical purposes. 
%\end{rmk}

\bigskip

\subsection{Definition of a global $\kappa$--entropy solution of \eqref{main}--\eqref{ini}} 
\begin{df}\label{Def1}
Let $\kappa$ be such that $0<\kappa<1$, the couple of functions $(\vr,\vu)$ is called a global $\kappa$--entropy solution to system \eqref{main}--\eqref{ini} if the following properties are satisfied:

\medskip

\noindent -- The mass equation is satisfied in the following sense 
\begin{equation} \label{mass}
- \intO {\vr \, \partial_t\xi} - \intO {\vr\vu\cdot \Grad \xi }= \intO {\vr^0\xi(0)}
\end{equation}
for all $\xi \in C^\infty_c([0,T)\times\Omega)$. 

\medskip

\noindent -- The momentum equation is satisfied in the following sense
\begin{equation}\label{u}
\begin{split}
&-\intO{\vr\vu\cdot\pt{\bf{\vcg{\phi}}}}
-\intO{(\vr \vu \otimes \vu):\Grad \vcg{\phi}}
+\intO{2 \mu(\vr)D(\vu):\Grad \vcg{\phi}}\\
&\quad +\intO{\lambda(\vr)\Div \, \vu\ \Div \, \vcg{\phi}}
-\intO{p(\vr)\Div \, \vcg{\phi}}=
\intO{\vr^0\vu^0\cdot\vcg{\phi}(0)}
\end{split}
\end{equation}
for all ${\bf{\vcg{\phi}}} \in (C^\infty_c([0,T)\times\Omega))^3$. 

\medskip

\noindent -- Moreover $(\vr,\vu)$ satisfies, for all $t\in [0,T]$,  the following $\kappa$-entropy estimates
\eq{\label{KappaEntropy}
& \sup_{t\in [0,T]} \Bigl[\intO{\vr\lr{\frac{|\vu+2\kappa\Grad\varphi(\vr)|^2}{2}
+ ( 1-\kappa)\kappa\frac{|2 \Grad\varphi(\vr)|^2}{2}  }(t)}
 +   \intO{ \vr e(\vr)(t)}\Bigr] \\
& +  2\kappa\int_0^T\intO{\mu(\vr) |A(\vu)|^2}\, {\rm d}s
 + 2 \kappa \int_0^T\intO {\frac{\mu'(\vr)p'(\vr)}{\vr} |\Grad\vr|^2}\, {\rm d}s  \\
& + 2(1-\kappa)\int_0^T\Bigl[\intO{\mu(\vr) |D(\vu)|^2} 
    + \intO{(\mu'(\vr)\vr-\mu(\vr))|\Div  \, \vu|^2}\Bigr] \, {\rm d}s\\
&\hskip2cm  \le \intO{\vr\lr{\frac{|\vu +2\kappa\Grad\varphi(\vr)|^2}{2}
+(1-\kappa)\kappa\frac{|2 \Grad\varphi(\vr)|^2}{2}  }(0)}
 +   \intO{ \vr_0 e(\vr_0)}
}
with $\varphi'(s)=\mu'(s)/s$,  $A(\vu)=\frac{1}{2}(\Grad\vu-\Grad^t\vu)$ and  the internal energy $e(\vr)$ defined by
$$\frac{\vr^2{\rm d} e(\vr)}{{\rm d}\vr}=p(\vr).$$
\end{df}

%%%I think it is well explained in the introduction, no need to repeat it%%%%%
%Let us note that such $\kappa$-entropy is a nonlinear generalization of the
%pseudocoercivity properties shown on linearized compressible Navier-Stokes
%equations. For reader's convience, we provide the calculations in this 
%framework in Section \ref{hypoc} after the proof of existence. 

\begin{rmk}
Note  that 
\eq{ \label{Enkappa}
\intO{\vr\lr{\frac{|\vu+2\kappa\Grad\varphi(\vr)|^2}{2}
+ (1-\kappa)\kappa\frac{|2\Grad\varphi(\vr)|^2}{2}  }(t)} \\=
 \intO{\vr\lr{(1-\kappa)\frac{|\vu|^2}{2}
+ \kappa\frac{|\vu + 2\Grad\varphi(\vr)|^2}{2}  }(t)}.
}
That means that the $\kappa$-entropy demonstrates  two-velocity structure in
the usual compressible Navier-Stokes system. For an introduction to the
two-velocity hydrodynamics and thermodynamics we refer to   {\rm  \cite{Sh}} and {\rm \cite{GaSh}} .  Indeed, using the identity
$$\vu + 2 \kappa\nabla\varphi(\vr) = (1-\kappa)\vu + \kappa (\vu+2\nabla\varphi(\vr))$$
we see that \eqref{Enkappa} is the kinetic energy of a two-fluid mixture having the velocities
$\vu$ and $\vu+2\nabla\varphi(\vr)$ and $\kappa$ playing the role of the mass fraction. 
 %The interested reader is referred to {\rm \cite{Sh}} for similar calculations.
%Before, the two first authors were using the energy estimate and the BD entropy  separately \cite{BrDe4, BrDe2}.  
\end{rmk}

\begin{rmk}
Equality \eqref{KappaEntropy} is a generalization of the BD entropy obtained by {\sc D. Bresch} and {\sc B. Desjardins} in the case $\kappa=1$. More precisely, they formally derived the following identity
\begin{equation} \label{BDentropy}
\Dt\intO{\vr \frac{|\vu + 2\Grad\varphi(\vr)|^2}{2} }
 +  \Dt \intO{ \vr e(\vr)}  
+ \intO{2 \mu(\vr) |A(\vw)|^2}
+ \intO {2 \frac{\mu'(\vr)p'(\vr)}{\vr} |\Grad\vr|^2} =0.
\end{equation}
If a global solution satisfies the BD entropy and the standard energy balance
$$ 
  \Dt\intO{\vr \frac{|\vu|^2}{2} }
 +  \Dt \intO{ \pi(\vr)}  
+ \intO{2 \mu(\vr) |D(\vu)|^2}  + \intO{\lambda(\vr) |\Div \, \vu|^2}=0,
$$
it also satisfies the $\kappa$-entropy estimate for all $0<\kappa<1$.
It suffices to use the identity
 $$ \kappa |\vu+2\Grad\varphi(\vr)|^2  + (1-\kappa)|\vu|^2 
    = |\vu+2\kappa\Grad\varphi(\vr)|^2 + (1-\kappa)\kappa |2\Grad\varphi(\vr)|^2
 $$
 and therfore to add $\kappa$ times the BD entropy to $(1-\kappa)$ times
 the energy.  A global weak solution satisfying the BD entropy is therefore a $\kappa$-entropy solution for all $0<\kappa<1$. The converse is however not clear.
 % because there is no reason for a $\kappa$-entropy solution 
%to be a global weak solution which satisfies the BD entropy.
\end{rmk}

%\begin{rmk} Note that, in three dimension,  we have the well-known relation 
% $$|D(\vu)|^2
%     = \left|D(\vu) - \frac{1}{3}\Div  \vu\,  \vc{I}\right|^2
%     + \frac{1}{3} |\Div  \vu|^2.$$
%This explains the assumption on $3\lambda(\vr) + 2 \mu(\vr)$ because
%$\lambda(\vr)  = 2 (\mu(\vr)\vr - \mu(\vr))$.
%\end{rmk}
\subsection{Main results} 

The initial data \eqref{ini} are assumed to satisfy
\eq{\vr^0\geq0,\quad \vr^0\in L^1(\Omega),\quad \vr^0e(\vr^0)\in L^1(\Omega),\quad \frac{\Grad\mu(\vr^0)}{\sqrt{\vr^0}}\in L^2(\Omega),\label{ini1}}
\eq{\frac{|\vc{m}^0|^2}{\vr^0}=0,\quad \text{a.e. on } \{x\in\Omega:\vr^0(x)=0\},\quad \frac{\vc{m}^0}{\sqrt{\vr^0}}\in L^2(\Omega)\label{ini2}}
and we consider the periodic boundary conditions
$$\Omega=\mathbb{T}^3.$$

\subsubsection{ Compressible Navier--Stokes equations with singular pressure}
In what follows we will make some assumptions concerning the viscosity coefficients and the pressure.
We assume that $\mu(\cdot)$, $\lambda(\cdot)$ are $C^1([0,\infty))$ such that $\mu'(\vr)\geq c>0$, 
$\mu(0)=0$, the following relation is satisfied
$$\lambda(\vr) = 2(\mu'(\vr)\vr - \mu(\vr)).$$
Moreover, there exists positive constants $c_0,c_1$, $\vr^*$, 
\eq{m>3/4, \quad 2/3<n<\frac{\gamma^-+1}{2},\label{mn}}
such that:
$$\text{for all}\quad s<\vr^*,\quad \mu(s)\geq c_0 s^n \quad \text{and} \quad3\lambda(s)+2\mu(s)\geq s^n,$$
$$\text{for all}\quad s\geq\vr^*,\quad c_1 s^m\leq\mu(s)\leq\frac{s^m}{c_1}\quad \text{and} \quad c_1 s^m\leq 3\lambda(s)+2\mu(s)\leq \frac{s^m}{c_1}.$$
The pressure \eqref{coldp} is a  $C^1([0,\infty))$ of $\vr$ and we assume that:
\eq{\gamma^+>1,\quad \gamma^->\frac{2n(3m-2)}{4m-3}-1.\label{gpm}}
Our main results reads as follows:
\begin{thm}\label{T_main} Assume that $0<\kappa < 1$ is fixed.
If the initial data, viscosity coefficients and $\gamma^{\pm}$ satisfy all the assumptions above
with a singular pressure given by \eqref{coldp} 
then there exists a global in time weak $\kappa$-entropy solution in the sense of Definition \ref{Def1}.
\end{thm}

\begin{rmk}
Assumptions \eqref{mn} and \eqref{gpm} are required to pass to the limit in the convective term.  More precisely one needs to ensure that
$$\vr^{\frac{1}{2}}\vu\in L^s(0,T; L^r(\Omega))$$
with $r,s>2$: the details of this estimate will be given in at the end of Section \ref{SS:Recov}.
The assumption  on $3\lambda(\vr) + 2 \mu(\vr)$ is used to get estimates from\eqref{KappaEntropy}. In the course of the proof, we combine BD relation $\lambda(\vr)  = 2 (\mu(\vr)\vr - \mu(\vr))$ with the following equality
 $$|D(\vu)|^2
     = \left|D(\vu) - \frac{1}{3}\Div  \vu\,  \vc{I}\right|^2
     + \frac{1}{3} |\Div  \vu|^2.$$
\end{rmk}

\subsubsection{ Compressible Navier--Stokes equations with a turbulent drag term}
 In this part, we assume the following hypotheses on viscosities $\mu(\vr)$ and $\lambda(\vr)$
 that may be found for instance in \cite{MeVa}:
let  $\mu(\cdot)$, $\lambda(\cdot)$ are $C^1([0,\infty))$ such that there exists a positive
$\nu\in(0,1)$ with 
\eq{\mu'(\vr)\geq \nu, \qquad \mu(0)\ge 0,\label{assum1}}
and
\eq{|\lambda'(\vr)| \le \frac{1}{\nu} \mu'(\vr), \qquad
    \nu\mu(\vr) \le 2\mu(\vr) + 3\lambda(\vr) \le \frac{1}{\nu}\mu(\vr).\label{assum2}}
The following relation (introduced by the two first authors) is also assumed 
\eq{\lambda(\vr) = 2(\mu'(\vr)\vr - \mu(\vr)).\label{rel}}
As stressed in \cite{MeVa},  the hypothesis above imply that 
\eq{C \vr^{2/3+\nu/3} \le \mu(\vr) \le C \vr^{2/3+1/(3\nu)} \hbox{ when } \vr\ge 1,\label{rhogrand}}
\eq{C\vr^{2/3+1/(3\nu)}\le \mu(\vr) \le C \vr^{2/3+\nu/3} \hbox{ when } \vr \le 1.\label{rhopetit}}
The interested reader is referred to \cite{MeVa} for more details to see where such hypothesis
are used for stability purposes.\\

\noindent In this part, we add the following hypothesis to deal with the turbulent drag term because,
compared to \cite{BrDeGe}, the $\kappa$-entropy mix the usual energy estimate and the BD-entropy
(less information is available):
\eq{ \mu(\vr) \le C \vr^{2/3+1/(3\eta)}, \hbox{ with } \eta \in (1/\gamma, 1) \hbox{ when } \vr \ge 1.
     \label{assum3}}

\begin{thm}\label{T_main2} Assume that $0<\kappa<1$ be fixed. Let us assume \eqref{assum1}--\eqref{assum3} and \eqref{ini1}-\eqref{ini2} be satisfied with a pressure $p(\vr)=a\vr^\gamma$
with $\gamma>1$, then there exists global in time  weak $\kappa$-entropy solution of \eqref{main2}
in the similar sense of Definition {\rm \ref{Def1}} (with the extra term coming from the turbulent drag and
the power-law pressure). 
\end{thm}

\begin{rmk} The drag term gives the extra information on $\vr|\vu|^2$ needed to
pass to the limit in the convective term without having the Mellet-Vasseur estimate. When no turbulent drag term is included  an extra assumption is required to get the 
Mellet-Vasseur estimate, namely, if $\gamma\ge3$, it is also assumed that
$$\liminf_{\vr\to +\infty} \frac{\mu(\vr)}{\vr^{\gamma/3+\zeta}} >0$$
with $\zeta>0$.
\end{rmk}

\begin{rmk}
The  system of compressible Navier-Stokes equations with BD structure with different boundary conditions was investigated in  {\rm \cite{BrDeGe}}. The authors considered  the Dirichlet condition for the momentum 
$$\vr \vu\vert_{\partial \Omega} = 0$$
as well as the Navier type condition 
$$\vr \vu\cdot \vc{n}\vert_{\partial\Omega} = 0, \qquad
   \mu(\vr) (D(\vu)\vc{n})_{\vcg{\tau}}\vert_{\partial \Omega} 
      =  - \alpha \mu(\vr) \vu_{\vcg{\tau}}\vert_{\partial \Omega},$$
together with additional boundary condition on the density, namely
 $$\mu(\vr) \nabla \varphi(\vr) \times \vc{n}\vert_{\partial\Omega} = 0.$$
 This type of boundary conditions could perhaps give transmission boundary conditions between the two quantities $\vu+2\kappa \nabla \varphi(\vr)$  and $2\nabla\varphi(\vr)$ and helps to conclude
 in bounded domains. 
\end{rmk}

%%%%%this has been explained in the introduction%%%%%
%\begin{rmk} Note that our constructive scheme is relevant for
%numerical purposes. It has  already been used in the inviscid framework by {\sc P. Noble} and
%{\sc J.-P.  Vila} in {\rm \cite{NoVi}} with extra dispersive terms (see also references cited therein).
%\end{rmk}

\subsection{Change of variable and $\kappa$--entropy}

   We generalize to the compressible framework some ideas developped recently in
\cite{BrGiZa} for low Mach number systems with heat conductivity effects. 
More precisely, let use define the following 
velocity field generalizing the one introduced  in the BD entropy estimate
\eq{\vw=\vu+ 2 \kappa \Grad\vp(\vr)} 
 with 
 \eq{\vp'(\vr) = \frac{\mu'(\vr)}{\vr}.\label{vp_mu}} 
  Note that BD entropy, as introduced in \cite{BrDe1} by the first two authors, corresponds to $\kappa=1$. 
  Assuming that solutions to \eqref{main} are smooth enough, it can be shown that $\vw$ satisfies the following evolution equation
\eq{
            & \pt\lr{\vr\vw}+\Div\lr{\vr\vu\otimes\vw}
            - 2 (1-\kappa)\Div (\mu(\vr)D(\vw)) - 2 \kappa \Div (\mu(\vr)A(\vw)) \\
           & +  4(1-\kappa)\kappa \Div (\mu(\vr)\Grad\Grad \varphi(\vr))   
            - \Grad (\bigl(\lambda(\vr)- 2\kappa (\mu'(\vr)\vr
            -\mu(\vr)\bigr)\Div \, \vu) 
            +\Grad P(\vr)=\vc{0}.
}        
 Let us now write the equation satisfied by $(\vr,\vw,\Grad\varphi(\vr))$.  We get the system
        \eq{\label{main2}  
     & \pt\vr + \Div\lr{\vr\vw} - 2 \kappa \lap \mu(\vr)= 0,\\
      &\pt\lr{\vr\vw} + \Div\lr{\vr\vu\otimes \vw}  -  2(1-\kappa) \Div\lr{\mu(\vr) \Grad \vw} 
          - 2 \kappa \Div(\mu(\vr) A(\vw)) \\ 
          &\qquad\qquad+ 4(1-\kappa)\kappa \Div(\mu(\vr)\Grad^2\varphi(\vr))
         -  \Grad (\bigl(\lambda(\vr)- 2\kappa (\mu'(\vr)\vr-\mu(\vr))\bigr)\Div \, \vu)+ \Grad P(\vr) = \vc{0},\\
     & \pt(\vr \Grad\varphi(\vr)) + \Div (\vr \vu\otimes \Grad\varphi(\vr)) 
      - 2 \kappa  \Div (\mu(\vr)\Grad^2\varphi(\vr))
     + \Div (\mu(\vr)\Grad^{t} \vw)  \\
      &\qquad\qquad + \Grad ((\mu'(\vr)\vr-\mu(\vr))\Div \, \vu) = \vc{0},       \\
         &\vw=\vu + 2\kappa\Grad \vp(\vr).
    } 
Taking the scalar product of the equation satisfied by $\vw$  with $\vw$,
the scalar product of the equation satisfied by $\Grad\varphi(\vr)$ with $4(1-\kappa)\kappa\Grad\varphi(\vr)$ and adding the resulting expressions  we get the $\kappa$--entropy 
\eq{\label{aa}
& \Dt\intO{\vr\lr{\frac{|\vw|^2}{2}
+ (1-\kappa)\kappa\frac{|2\Grad\varphi(\vr)|^2}{2}  }}
 +  \Dt \intO{ \vr e(\vr)} \\
& +  2\kappa\intO{\mu(\vr) |A(\vw)|^2}
 + 2\kappa \intO {\frac{\mu'(\vr)p'(\vr)}{\vr} |\Grad\vr|^2}  \\
& + 2(1-\kappa)\Bigl[\intO{\mu(\vr) |D(\vu)|^2} 
    + \intO{(\mu'(\vr)\vr-\mu(\vr))|\Div \, \vu |^2}\Bigr]
=0
}
where we used assumption \eqref{BD} to write
\eq{
& \intO{(\lambda(\vr)-2\kappa(\mu'(\vr)\vr-\mu(\vr)) \Div \, \vu \,  \Div \, \vw
       + 4(1-\kappa)\kappa (\mu'(\vr)\vr-\mu(\vr)) \Div \, \vu \, \, \Div \Grad\varphi(\vr)} \\
&       = 2(1-\kappa) \intO{(\mu'(\vr)\vr-\mu(\vr)) |\Div\,\vu|^2}
}
recalling that $\vu = \vw - 2\kappa\Grad\varphi(\vr)$.

\begin{rmk}
 Relation \eqref{aa}  generalizes the one obtained in  \cite{Ju} (see also \cite{GiVi}) for quantum Navier--Stokes equations, it is enough to take $\kappa=1/2$ and $\varphi=\log\vr$. However, in contrast to aforementioned papers, we do not need the  Bohm potential formula to conclude.
% Our definition of global in time $\kappa$--entropy solutions will be linked to \eqref{aa} --
% A solution which satisfies the mass and momentum equations in a weak sense and the
% $\kappa$-entropy estimate. %%%%%%The definition is already stated
\end{rmk}

\section{Construction of solution}\label{S:constr}

    Following the idea developped recently by {\sc D. Bresch, V. Giovangigli} and {\sc E. Zatorska}   in \cite{BrGiZa} for the following low Mach system
     \begin{equation}\label{Limref2}
\begin{array}{c}
\partial_t \vr + \Div (\vr \vu) = 0, \\
\pt (\vr \vu) + \Div (\vr \vu \otimes \vu)
         +\Grad \pi =  2 \Div (\mu(\vr) D(\vu)) + \Grad (\lambda(\vr) \Div \, \vu) ,\\
\Div \, \vu = -  2 \kappa \Delta \varphi(\vr),
\end{array}
\end{equation} 
with $\varphi$ an increasing function of $\vr$ and $0<\kappa<1$ fixed. We construct the approximate solution to system \eqref{main2} using an augmented approximate system. More precisely, we introduce a new unknown $\vv$, which is not yet known to satisfy $\vv=2\Grad\vp(\vr)$. Our aim will be to find a solution $(\vr,\vw,\vv)$ of the following system:
\eq{\label{approx}
&\pt\vr+ \Div(\vr \vw) - 2 \kappa\lap\mu( \vr) = 0,\\
& \pt\lr{\vr \vw}+ \Div  \lr{(\vr \vw - 2 \kappa  \Grad\mu( \vr)) \otimes \vw}  -  \Grad \lr{\lr{\lambda(\vr)- 2\kappa (\mu'(\vr)\vr-\mu(\vr))}\Div(\vw- \kappa\vv)} \\
 &- 2(1-\kappa)  \Div (\mu(\vr) D(\vw)) - 2 \kappa \Div (\mu(\vr) A(\vw))  
 + \Grad p(\vr)= - 2\kappa (1-\kappa) \Div(\mu(\vr) \Grad\vv),\\
&\pt(\vr \vv) + \Div((\vr \vw- 2 \kappa \Grad \mu( \vr))\otimes \vv) 
-2 \kappa \Div(\mu(\vr)\Grad\vv) + 2 \Grad ((\mu'(\vr)\vr - \mu(\vr))\Div \lr{\vw- \kappa\vv})\\
&= - 2 \Div(\mu(\vr) \Grad^{t} \vw),
} 

%%%%This is explained in the introduction
 \begin{rmk} To see the link between the above system and the one from {\rm \cite{Sh}}, one should take the main velocity equal to $\vw= \vu+2\kappa\Grad \vp(\vr)$ and the drift equal to $\vv= 2\Grad \varphi(\vr)$.
 \end{rmk}
 
 \medskip
 
\noindent Our idea of construction of solution is based on the following $\kappa$-entropy equality
    \eq{\label{aa-regu}
&\Dt\intO{\vr\lr{\frac{|\vw|^2}{2}
+ (1-\kappa)\kappa\frac{|\vv|^2}{2}  }}
 +  \Dt \intO{ \vr e(\vr)} + 2(1-\kappa)\intO{\mu(\vr) |D(\vw) -\kappa\Grad \vv|^2} \\
&\qquad+ 2\kappa\intO{\mu(\vr) |A(\vw)|^2}
+ 2(1-\kappa)\intO{(\mu'(\vr)\vr-\mu(\vr))|\Div (\vw-\kappa\vv)|^2}
\\
&\qquad+ 2 \intO {\frac{\kappa\mu'(\vr)p'(\vr)}{\vr} |\Grad\vr|^2} =0,
}
that holds for any sufficiently smooth solution of \eqref{approx}.
   In order to build such a  solution for $\varepsilon>0$ given, we need to go through several levels of approximations. For example, to build a solution of the nonlinear parabolic equation for $\vr$, some assumptions are required on the coefficients. Two smoothing parameters $\alpha>0$ and $\delta>0$(denoting standard mollification with respect to $t$ and $x$) are therefore introduced in all the transport terms, so that the  approximate system can be rewritten as
\eq{\label{approx2}
&\pt\vr+ \Div(\vr [\vw]_\delta) - 2 \kappa\Div\lr{[\mu'( \vr)]_\alpha \Grad\vr} = 0,\\
 &\pt\lr{\vr \vw}+ \Div  ((\vr [\vw]_\delta - 2 \kappa  [\mu'( \vr)]_\alpha \Grad\vr) \otimes \vw)  -  \Grad (\lr{\lambda(\vr)- 2 \kappa (\mu'(\vr)\vr-\mu(\vr))}\Div(\vw- \kappa\vv)) \\
&\qquad\qquad- 2(1-\kappa)  \Div (\mu(\vr) D(\vw)) - 2 \kappa \Div (\mu(\vr) A(\vw))  
 + \varepsilon \lap^{2s} \vw - \varepsilon \Div ((1+|\Grad\vw|^2)\Grad\vw)  + \Grad p(\vr)\\
&\qquad\quad = -  2 \kappa (1-\kappa) \Div(\mu(\vr) \Grad\vv),\\
&\pt(\vr \vv) + \Div((\vr [\vw]_\delta - 2 \kappa [\mu'( \vr)]_\alpha \Grad\vr)\otimes \vv) \\
&\qquad\qquad -2 \kappa \Div(\mu(\vr)\Grad\vv) + 2\Grad \lr{(\mu'(\vr)\vr - \mu(\vr))\Div \lr{\vw-\kappa\vv}}\\
&\qquad\quad= - 2\Div(\mu(\vr) \Grad^{t} \vw).
}
 Compared to \cite{BrGiZa}, we have an extra  term in the continuity equation because $\vw$ is no longer divergence free. In addition, a smoothing high-order derivative term $\lap^{2s}\vw$ with $s\geq 2$, depending on small parameter $\ep>0$ has to be introduced to control large spatial variations of $\vw$,  because $\Div \, \vw$ is not {\it a-priori}  bounded in $L^1(0,T;(L^\infty(\Omega))^3$. 
    Such bound will be required to be able to have bounds on the density. 
    We will also need to show that $\vv =2 \Grad \varphi(\vr)$ at some point of the construction process
and the second term in the regularization process will be helpful.

In order to solve \eqref{approx2} for given $\alpha,\delta, \ep>0$, the equations for the two velocities $(\vw,\vv)$ are projected onto finite dimensional spaces $(X_n,Y_n)$ as is classically done when using Faedo-Galerkin approximation.
   The first step is to prove global in time existence of solutions when $n$, $\alpha$, $\delta$ and $\varepsilon$ are fixed, then  to pass to the limit $\alpha\to 0$ and after to let  $n\to\infty$. Next, we let $\delta\to 0$ in order to prove that $\vv = 2 \Grad \varphi(\vr)$.  At the end, we combine the equation satisfied by $\vw$ and $\vv$ to get an equation satisified by $\vu$ (multiply the equation satisfied by $\vv$ by $\kappa$ and subtract from the equation for $\vw$). 
The last limit passage $\varepsilon\to0$ is performed  already for the system written in terms of  $\vu$ which gives the weak formulation of the original system \eqref{main}. 
This last step follows the same lines as the proof presented by the two first authors for the heat-conducting case, see \cite{BrDe3}. It essentially uses the presence of singular pressure, which allows to stabilize the system close to vacuum or the presence of the turbulent drag term to have enough control on the velocity field. More complicated construction in the case of Navier-Stokes type model for compressible mixture with viscosity $\mu(\vr)=\vr$ can be found in the PhD thesis of the third author \cite{Za}, see also \cite{MuPoZa1, MuPoZa2}.

 \medskip
 
 \subsection{Existence of solutions for the full approximation} Below we present the basic level of approximation procedure.

1. The continuity equation is replaced by its regularized version
\eq{&\pt\vr+ \Div(\vr [\vw]_\delta) - \kappa\, \Div([\mu'({\vr})]_\alpha \Grad \vr) = 0,\\
&\qquad\qquad
\vr(0,x)=[\vr^0]_{\delta},
\label{A_cont}}
where $\alpha,\delta$ denote the standard regularizations with respect to time and space.
%With such smoothing properties, we can use the estimates given by Theorem 1 which are uniform with 
%respect to $\delta$ when far from vacuum.  This will be used at different stages: the construction and
%the proof that $\vv = 2\Grad \varphi(\vr)$ before passing to the limit with respect to $\varepsilon$.
Note that similar double regularization have been recently used in \cite{Li} and discussed in the book by J.L. {\sc Vazquez} \cite{Va}.

\smallskip

2. The momentum equation is replaced by its Faedo-Galerkin approximation with additional regularizing term $\ep\bigl[\Delta^{2s} \vw - \Div((1+|\Grad\vw|^2)\Grad\vw)]$
\eq{\label{FG}
&\intO{\vr\vw(\tau)\cdot\vcg{\phi}}
-\inttauO{((\vr [\vw]_\delta - 2\kappa[\mu'({\vr})]_\alpha \Grad \vr) \otimes \vw):\Grad \vcg{\phi}}\\
&+ 2(1-\kappa)\inttauO{\mu(\vr)D(\vw):\Grad \vcg{\phi}} + 2\kappa\inttauO{\mu(\vr)A(\vw): \Grad\vcg{\phi}} \\
&+ 2(1-\kappa)\inttauO{\lr{\mu'(\vr)\vr -\mu(\vr) }\Div \, \vw\, \Div \, \vcg{\phi}} \\
& - 2 \kappa(1-\kappa)\inttauO{\mu(\vr)\Grad\vv:\Grad\vcg{\phi}}  -  2\kappa(1-\kappa) \inttauO{(\mu'(\vr)\vr -\mu(\vr) )\Div \, \vv\, \Div \, \vcg{\phi}}  \\
&-\inttauO{p(\vr)\Div \, \vcg{\phi}}+\ep\inttauO{\lr{\lap^s\vw\cdot\lap^s\vcg{\phi}+(1+|\Grad\vw|^2) \Grad\vw:\Grad\vcg{\phi}}}
\\
&=\intO{(\vr\vw)^{0}\cdot\vcg{\phi}},
}
satisfied for any $\tau\in[0,T]$ and any test function $\vcg{\phi}\in X_{n}$, where 
$X_{n}=\operatorname{span}\{\vcg{\phi}_{i}\}_{i=1}^{n}$  and $\{\vcg{\phi}_{i}\}_{i=1}^{\infty}$ is an orthonormal basis in $(W^{1,2}(\Omega))^3$ with 
$\vcg{\phi}_{i}\in ({\cal C}^\infty(\Omega))^3$ for all $i\in N$.

\smallskip

3. The Faedo-Galerkin approximation for the artificial equation
\eq{\label{FGtheta}
&\intO{\vr\vv(\tau)\cdot {\bf{\vcg{\xi}}}}
-\inttauO{((\vr [\vw]_\delta - 2\kappa [\mu'({\vr})]_\alpha \Grad \vr) \otimes \vv):\Grad \vcg{\xi}}
\\
&+2 \kappa\inttauO{\mu(\vr)\Grad\vv:\Grad \vcg{\xi}}
+2 \kappa\inttauO{(\mu'(\vr)\vr - \mu(\vr))\Div \, \vv\,\Div \, \vcg{\xi}}
\\
&-2 \inttauO{(\mu'(\vr)\vr - \mu(\vr))\Div \, \vw\,\Div \, \vcg{\xi}}
-2\inttauO{\mu(\vr)\Grad^t\vw:\Grad\vcg{\xi}}\\
&=\intO{(\vr\vv)^{0} \cdot \vcg{\xi}},
}
satisfied for any $\tau\in[0,T]$ and any test function $\vcg{\xi}\in Y_{n}$, where $Y_{n}=\operatorname{span}\{\vcg{\xi}_{i}\}_{i=1}^{n}$   and $\{\vcg{\xi}_{i}\}_{i=1}^{\infty}$ is an orthonormal basis in $(W^{1,2}(\Omega))^3$
with  $\vcg{\xi}_{i}\in ({\cal C}^\infty(\Omega))^3$ for all $i\in N$.

\bigskip

\noindent{\bf Existence of solutions to the continuity equation.} For fixed $\vw\in C([0,T]; X_n)$  we solve the continuity equation, which is now quasi-linear parabolic equation with smooth coefficients.
Thus, application of classical existence theory of Lady{\v{z}}enskaja, Solonnikov and Uralceva \cite{LSU} (see for example Theorem 10.24 from \cite{FeNo}, which is a combination of Theorems 7.2, 7.3 and 7.4 from \cite{LSU}) yields the following result:
\begin{thm}\label{LSU} 
Let $\nu\in (0,1)$ and suppose that the initial condition $\vr_\delta^0\in C^{2+\nu}(\Ov\Omega)$ is such that $0<r\leq \vr_\delta^0\leq R$ and it satisfies the periodic boundary conditions. Then  problem \eqref{A_cont} possesses a unique classical solution $\vr$ from the class 
\begin{equation}\label{regvr}
V_{[0,T]}=\left\{\begin{array}{rl}
\vr&\in C([0,T];C^{2+\nu}(\Omega))\cap C^1([0,T]\times\Omega),\\
\pt\vr&\in C^{\nu/2}\left([0,T];C(\Omega)\right)
\end{array}\right\}
\end{equation}
and satisfying classical maximum principle 
\eq{0< r\leq\vr(t,x)\leq R. \label{max_gal}}
Moreover, the mapping $\vw\mapsto\vr(\vw)$ maps bounded sets in $C([0,T];X_{n})$ into bounded sets in $V_{[0,T]}$ and is continuous  with values in $C\big([0,T];C^{2+\nu'}(\Omega)\big)$, $0<\nu'<\nu<1$.
\end{thm}
\bigskip

\noindent{\bf Local existence of solutions to the Galerkin approximations.} Here we proceed as in the analogous proof performed in \cite{BrGiZa} for $\vw$ being the divergence-free vector field. In what follows, we will show that the integral equations \eqref{FG} and \eqref{FGtheta} possess the unique solution on possibly short time interval via fixed point argument.
More precisely, we will prove that there exists  time $T=T(n)$ and $(\vw,\vv)\in C([0,T];X_{n})\times C([0,T]; Y_n)$ satisfying (\ref{FG}, \ref{FGtheta}).
To this purpose let us rewrite this equations as a fixed point problem
\eq{\label{T}
(\vw(t),\vv(t))&=\lr{ {\mathcal{M}}_{\vr(t)}\left[P_{X_n}(\vr\vw)^0+\int_{0}^{t}{{\mathcal{K}}(\vw)(s) {\rm d}s}\right],\ \mathcal{N}_{\vr(t)}\left[P_{Y_n}(\vr\vv)^0+\int_{0}^{t}{{\mathcal{L}}(\vv)(s) {\rm d}s}\right]}\\
&={\mathcal{T}}[\vw,\vv](t)
}
where $\vr=\vr(\vw)$ is a solution to the continuity equation with $\vw$ given, $X_n^*$ is identified with $X_n$, so the symbol $\langle\cdot,\cdot\rangle_{(X_n,X_n)}$ denotes the action of a functional from $X_n^*(=X_n)$ on the element from $X_n$, similarly for $\langle\cdot,\cdot\rangle_{(Y_n,Y_n)}$, and 
$${\cal{M}}_{\vr(t)}:X_{n}\rightarrow X_{n},\quad\intO{\vr{\cal{M}}_{\vr(t)}[\vcg{\phi}]\cdot\vcg{\psi}}=\langle\vcg{\phi},\vcg{\psi}\rangle_{(X_n,X_n)}, \quad\vcg{\phi},\vcg{\psi}\in X_{n},$$
$${\cal{N}}_{\vr(t)}:Y_{n}\rightarrow Y_{n},\quad\intO{\vr{\cal{N}}_{\vr(t)}[\vcg{\xi}]\cdot\vcg{\zeta}}=\langle\vcg{\xi},\vcg{\zeta}\rangle_{(Y_n,Y_n)}, \quad\vcg{\xi},\vcg{\zeta}\in Y_{n},$$
$P_{X_n}$, $P_{Y_n}$ denote the projections of $L^2(\Omega)$ onto $X_n$, $Y_n$, respectively, and ${\cal K}(\vw), {\cal L}(\vv)$ are two operators defined as follows
\eqh{
{\cal K}: &\quad X_n\to X_n,\\
\langle{\cal K}(\vw),\vcg{\phi}\rangle_{(X_n,X_n)}=&\intO{((\vr [\vw]_\delta - 2\kappa [\mu'({\vr})]_\alpha \Grad \vr) \otimes \vw):\Grad \vcg{\phi}}
-2(1-\kappa)\intO{\mu(\vr)D(\vw):\Grad \vcg{\phi}}\\
&-2(1-\kappa)\intO{\lr{\mu'(\vr)\vr -\mu(\vr) }\Div  \, \vw\, \Div \, \vcg{\phi}} \\
&-2\kappa\intO{\mu(\vr)A(\vw): \Grad\vcg{\phi}}+2\kappa(1-\kappa)\intO{\mu(\vr)\Grad\vv:\Grad\vcg{\phi}}\\
&+ 2\kappa(1-\kappa) \intO{(\mu'(\vr)\vr -\mu(\vr) )\Div \, \vv\, \Div \, \vcg{\phi}} +\intO{p(\vr)\Div \,\vcg{\phi}}\\
&-\ep\intOB{\lap^s\vw\cdot\lap^s\vcg{\phi}+(1+|\Grad\vw|^2) \Grad\vw:\Grad\vcg{\phi}},
}
\eqh{
{\cal L}: &\quad Y_n\to Y_n,\\
\langle{\cal L}(\vv),\vcg{\xi}\rangle_{(Y_n,Y_n)}=
&\intO{((\vr [\vw]_\delta - 2\kappa [\mu'({\vr})]_\alpha \Grad \vr) \otimes \vv):\Grad \vcg{\xi}}
-2\kappa\intO{\mu(\vr)\Grad\vv:\Grad \vcg{\xi}}\\
&-2\kappa\intO{(\mu'(\vr)\vr - \mu(\vr))\Div \, \vv\ \Div \, \vcg{\xi}}
+2 \intO{(\mu'(\vr)\vr - \mu(\vr))\Div \, \vw\,\Div \, \vcg{\xi}}\\
&+2\intO{\mu(\vr)\Grad^t\vw:\Grad\vcg{\xi}}.}
Since $\vr(t,x)$ is bounded from below by a positive constant, we have
\eqh{\|{\cal M}_{\vr(t)}\|_{L(X_n,X_n)},\ \|{\cal N}_{\vr(t)}\|_{L(Y_n,Y_n)}\leq\frac{1}{r}.}
Moreover
\eq{\|{\cal M}_{\vr^1(t)}-{\cal M}_{\vr^2(t)}\|_{L(X_n,X_n)}
+ \|{\cal N}_{\vr^1(t)}-{\cal N}_{\vr^1(t)}\|_{L(Y_n,Y_n)}\leq c(n,r^1,r^2)\|\vr^1-\vr^2\|_{L^1(\Omega)},
\label{XY}
}
 and by the equivalence of norms on the finite dimensional space we prove that
 \eq{\|{\cal K}(\vw)\|_{X_n}+ \|{\cal L}(\vv)\|_{Y_n}\leq c(r,R,\|\Grad\vr\|_{L^2(\Omega)},\|\vw\|_{X_n}, \|\vv\|_{Y_n} ).\label{KL}}
Next, we consider a ball ${\cal B}$ in the space $C([0,\tau];X_{n})\times C([0,\tau];Y_{n})$:
	$${\cal B}_{M,\tau}=\left\{(\vw,\vv)\in C([0,\tau];X_{n})\times C([0,\tau];Y_{n}):\|\vw\|_{C([0,\tau];X_n)}+\|\vv\|_{C([0,\tau];Y_n)}\leq M\right\}.$$
Using estimates \eqref{XY}, \eqref{KL}, \eqref{regvr} and \eqref{max_gal}, one can check that ${\cal T}$ is a continuous mapping of the ball ${\cal B}_{M,\tau}$ into itself and for sufficiently small $\tau=T(n)$ it is a contraction. Therefore, it possesses a unique fixed point which is a solution to \eqref{FG} and \eqref{FGtheta} for $T=T(n)$.
\bigskip

\noindent{\bf Global existence of solutions.} In order to extend the local in-time solution obtained above to the global in time one, we need to find uniform (in time) estimates, so that the above procedure can be iterated.
First let us note, that $\vw,\vv$ obtained in the previous paragraph have better regularity with respect to time. It follows by taking the time derivative of \eqref{T} and using the estimates \eqref{regvr}, \eqref{max_gal}, that
$$(\vw,\vv)\in C^1([0,\tau];X_{n})\times C^1([0,\tau];Y_{n}).$$
This is an important feature since now we can take time derivatives of \eqref{FG} and \eqref{FGtheta} and use the test functions $\vcg{\phi}=\vw$ and $\vcg{\xi}=\vv$, respectively. We then obtain 
\begin{equation}\label{w1}
\begin{split}
&\Dt\intO{\vr\frac{|\vw|^2}{2}}+  \Dt \intO{ \vr e(\vr)}
+2(1-\kappa)\intO{\mu(\vr)|D(\vw)|^2}\\
&+2\kappa\intO{\mu(\vr)|A(\vw)|^2}+ 2(1-\kappa)\intO{\lr{\mu'(\vr)\vr -\mu(\vr) }(\Div \, \vw)^2} \\
&-2\kappa(1-\kappa)\intO{\mu(\vr)\Grad\vv:\Grad\vw}-  2\kappa(1-\kappa) \intO{(\mu'(\vr)\vr -\mu(\vr) )\Div \, \vv\, \Div \, \vw} \\
&+2\kappa\intTO {\frac{\mu'(\vr)p'(\vr)}{\vr} |\Grad\vr|^2} +\intO{p(\vr)\Div\lr{[\vw]_\delta-\vw}}\\
&+\ep\intOB{|\lap^s\vw|^2+(1+|\Grad\vw|^2) |\Grad\vw|^2}=0,
\end{split}
\end{equation}
and
\eq{\label{t1}
&\Dt\intO{\vr\frac{|\vv|^2}{2}}
+2\kappa\intTO{\mu(\vr)|\Grad\vv|^2}+2\kappa\intTO{(\mu'(\vr)\vr - \mu(\vr))(\Div \, \vv)^2}\\
&\quad
-2 \intO{(\mu'(\vr)\vr - \mu(\vr))\Div \, \vw\,\Div \, \vv}
-2\intTO{\mu(\vr)\Grad^t\vw:\Grad\vv}=0.
}
Therefore, multiplying \eqref{t1} by $(1-\kappa)\kappa$ and adding it to \eqref{w1}, w obtain
\eq{\label{aa1}
\Dt&\intO{\vr\lr{\frac{|\vw|^2}{2}
+(1-\kappa)\kappa\frac{|\vv|^2}{2}}}+  \Dt \intO{ \vr e(\vr)}+2(1-\kappa)\intO{\mu(\vr) |D(\vw) -\kappa\Grad \vv|^2} \\
&+2(1-\kappa)\intO{\lr{\mu'(\vr)\vr -\mu(\vr) }(\Div \, \vw-\kappa\Div \, \vv)^2}+2\kappa\intO{\mu(\vr) |A(\vw)|^2}\\
&+ 2\kappa \intO {\frac{\mu'(\vr)p'(\vr)}{\vr} |\Grad\vr|^2} +\ep\intOB{|\lap^s\vw|^2+(1+|\Grad\vw|^2) |\Grad\vw|^2}\\
&=-\intO{p(\vr)\Div\lr{[\vw]_\delta-\vw}}.
}
Integrating the above estimate with respect to time, using the H\'older and the Gronwall inequalities,
we obtain uniform estimate for $\vw$ and $\vv$ necessary 
to repeat the procedure described in the previous paragraph. Thus,  we obtain a global in time unique solution $(\vr,\vw,\vv)$ satisfying equations (\ref{A_cont}, \ref{FG}, \ref{FGtheta}).

\bigskip
\noindent{\bf Uniform estimates.} Below we  present uniform estimates that will allow us to pass to the limit with $\alpha$ and $n$ respectively.
First observe that multiplying continuity equation \eqref{A_cont} by $\vr_\alpha$ and integrating by parts with respect to $x$ gives
$$\frac{1}{2} \Dt \intO {|\vr_\alpha|^2} + 2\kappa \intO {[\mu'(\vr_\alpha)]_\alpha |\Grad\vr_\alpha|^2} = - \intO {\Div  [\vw]_\delta |\vr_\alpha|^2}
\le \|\Div [\vw]_\delta\|_{L^\infty(\Omega)} \intO{|\vr_\alpha|^2}
$$
Integrating this equality with respect to time and using the uniform bound of $\Div [\vw]_\delta$ in 
$L^1(0,T;L^\infty(\Omega))$  provides the following estimates 
\eq{\|\vr_\alpha\|_{L^\infty(0,T; L^2(\Omega))}+\|\sqrt{[\mu'(\vr_\alpha)]_\alpha }\Grad\vr_\alpha\|_{L^2(0,T;L^2(\Omega))}\leq c.\label{rho1}}
Moreover, the standard maximum principle gives boundedness of $\vr_\alpha$ from above and below. Indeed, multiplying equation \eqref{A_cont} by
$$
  \vr_\alpha^-=\max(0,r-\vr_\alpha)\quad\text{and}\quad \vr_\alpha^+=\min(0,R-\vr_\alpha),
$$
respectively we obtain also using the bound on $[\vw]_\delta$ that 
\eq{0< r_\varepsilon \leq\vr_\alpha(t,x)\leq R_\varepsilon. \label{max}}
%\begin{rmk}
%To prove these bounds one needs to know that $\vr_\alpha\in L^2(0,T;W^{1,2}(\Omega))$, which doesn't follow from \eqref{rho1}. This problem could be solved by adding a small viscosity parameter $\alpha$ and considering $\widetilde{[\mu'(\vr)]_\alpha}=[\mu'(\vr)]_\alpha+\alpha$ in place of $[\mu'(\vr)]_\alpha$.
%\end{rmk}
Next, using \eqref{max} and integrating \eqref{aa1} with respect to time we see that for $0<\kappa<1$ we have 
\eq{&\|\vw_\alpha\|_{L^\infty(0,T;L^2(\Omega)}
+\|\vw_\alpha\|_{L^2(0,T;H^{2s}(\Omega))}+
\|\Grad\vw_\alpha\|_{L^4(0,T;L^4(\Omega))}\\
&+\|\vv_\alpha\|_{L^\infty(0,T;L^2(\Omega))}
+\|\vv_\alpha\|_{L^2(0,T;H^{1}(\Omega))}\leq c.
\label{ueu}}
 In the above estimate the constant $c$ is uniform with respect to all approximation parameters 
 except $\ep$.

\subsection{Passage to the limit with respect to $\alpha$ and
with respect to $n$.}

\medskip

\noindent{\bf {Passage to the limit $\alpha\to0$.}} On the finite dimensional subspace all the norms are equivalent, therefore the space compactness of $\vw_\alpha$ and $\vv_\alpha$ is automatic. In fact, for $n$ fixed we also know that $\pt\vw_\alpha$ is bounded in $L^2(0,T; X_n)$, thus up to the extraction of a subsequence, $\vw_\alpha\to\vw$ strongly in $L^2(0,T;X_n)$ and the same can be deduced for $\vv_\alpha$.
The biggest problem is thus to pass to the limit in the term
\eq{[\mu'(\vr_\alpha)]_\alpha\Grad\vr_\alpha\otimes\vw_\alpha\label{conv_w} \quad \hbox{when} \quad 
\alpha \to 0}
which requires the strong convergence of the density and the weak convergence of the gradient of density. Observe that  higher order estimates for $\vr_\alpha$ cannot be used here due to $\alpha$-regularization of the coefficient $\mu'(\vr_\alpha)$, however the information from \eqref{max} and \eqref{rho1} is enough. 
   Indeed, using \eqref{max} and the assumption that $\mu'(\vr)\geq c>0$, estimate \eqref{rho1} implies that that up to a subsequence
\eqh{\vr_\alpha\to\vr\quad\text{weakly in }L^2(0,T; W^{1,2}(\Omega)) \quad \hbox{when} \quad 
\alpha \to 0,}
moreover $\pt\vr_\alpha\in L^2(0,T; W^{-1,2}(\Omega))$ and $\vr_\alpha \in L^\infty((0,T)\times\Omega)$, thus the Aubin-Lions lemma implies strong convergence of $\vr_\alpha$
\eqh{\vr_\alpha\to\vr\quad\text{strongly in }L^p((0,T)\times\Omega) \quad \hbox{when} \quad 
\alpha \to 0, \quad p<\infty.}
This justifies the passage to the limit in \eqref{conv_w}. Therefore, one is able to pass to the limit $\alpha\to 0$ in both velocity equations \eqref{FG} and \eqref{FGtheta}. The limit functions $(\vr,\vw,\vv)=(\vr_n,\vw_n,\vv_n)$ satisfy the following system of equations:
\begin{itemize}
\item the momentum equation
\begin{equation}\label{FG_a}
\begin{split}
&\langle\pt\lr{\vr_n\vw_n}(t),\vcg{\phi}\rangle_{(X_n^*,X_n)}
-\intO{((\vr_n [\vw_n]_\delta - 2\kappa \Grad \mu(\vr_n)) \otimes \vw_n)(t):\Grad \vcg{\phi}}\\
&+2(1-\kappa)\intO{\mu(\vr_n)D(\vw_n)(t):\Grad \vcg{\phi}}
+2\kappa\intO{\mu(\vr_n)A(\vw_n)(t): \Grad\vcg{\phi}}\\
&+ 2(1-\kappa)\intO{\lr{\mu'(\vr_n)\vr_n -\mu(\vr_n) }\Div \, \vw_n(t)\, \Div \, \vcg{\phi}} \\
&-2\kappa(1-\kappa)\intO{\mu(\vr_n)\Grad\vv_n(t):\Grad\vcg{\phi}}
 -  2\kappa(1-\kappa) \intO{(\mu'(\vr_n)\vr_n -\mu(\vr_n) )\Div \, \vv_n(t)\, \Div \, \vcg{\phi}}  \\
&-\intO{p(\vr_n)(t)\Div \, \vcg{\phi}}+\ep\intOB{\lap^s\vw_n(t)\cdot\lap^s\vcg{\phi}+(1+|\Grad\vw_n|^2) \Grad\vw_n(t):\Grad\vcg{\phi}}=0
\end{split}
\end{equation}
satisfied for all $\vcg{\phi}\in X_n$ with $t\in[0,T]$, 
\item the auxiliary equation for $\vv_n$
\eq{\label{FGtheta_a}
&\langle\pt\lr{\vr_n\vv_n}(t),\vcg{\xi}\rangle_{(Y_n^*,Y_n)}
-\intO{((\vr_n [\vw_n]_\delta - 2\kappa  \Grad \mu(\vr_n)) \otimes \vv_n)(t):\Grad \vcg{\xi}}\\
&+2\kappa\intO{\mu(\vr_n)\Grad\vv_n(t):\Grad \vcg{\xi}}
+2\kappa\intO{(\mu'(\vr_n)\vr_n - \mu(\vr_n))\Div \, \vv_n(t)\, \Div \, \vcg{\xi}}\\
&-2 \intO{(\mu'(\vr_n)\vr_n - \mu(\vr_n))\Div \, \vw_n(t)\,\Div \, \vcg{\xi}}
-2\intO{\mu(\vr_n)\Grad^t\vw_n(t):\Grad\vcg{\xi}}=0
}
is satisfied for all $\vcg{\xi}\in Y_n$ with $t\in[0,T]$.
\end{itemize}
However, so far we only know that the approximate continuity equation \eqref{A_cont} is satisfied in the sense of distributions 
\eq{\label{weak1}
\intT{\langle\pt\vr_n,\phi\rangle_{(W^{-1,2}(\Omega), W^{1,2}(\Omega))}}-\intTO{\vr_n[\vw_n]_\delta\cdot\Grad\phi}+\kappa\intTO{\Grad\mu(\vr_n)\cdot\Grad\phi}=0}
for any test function $\phi$ from $L^2(0,T;W^{1,2}(\Omega))$.
But on the other hand, we know that  in the sense of distributions
$$
  \psi=\pt\vr_n-\kappa\Div\lr{\mu'(\vr_n)\Grad\vr_n}=-\vr_n\Div[\vw_n]_\delta-[\vw_n]_\delta\cdot\Grad\vr_n\in L^2((0,T)\times\Omega),$$
if  $\vw_n\in L^\infty(0,T; L^\infty(\Omega))$. Indeed, note that we have
\eq{\|\vr_n\Div[\vw_n]_\delta\|_{L^2(0,T; L^2(\Omega))}\leq\|\vr_n\|_{L^\infty(0,T; L^\infty(\Omega))}\|\Div[\vw_n]_\delta\|_{L^2(0,T; L^2(\Omega))}
}
\eq{\|[\vw_n]_\delta\cdot\Grad\vr_n\|_{L^2(0,T; L^2(\Omega))}\leq\|[\vw_n]_\delta\|_{L^\infty(0,T; L^\infty(\Omega))}\|\Grad\vr_n\|_{L^2(0,T; L^2(\Omega))}}
and the r.h.s. of above inequalities is bounded at the level of Galerkin approximations, therefore also
$$\psi\mu'(\vr_n)\in L^2((0,T)\times\Omega).$$
Taking the product of $\psi$ and $\psi\mu'(\vr_n)$ we obtain
\eqh{\intO{\mu'(\vr_n)\lr{\pt\vr_n-\kappa\lap\mu(\vr_n)}^2}\leq c(t)\in L^1((0,T))}
and the above integral gives rise to estimates
\eq{\intO{\mu'(\vr_n)(\pt\vr_n)^2}+\kappa^2\intO{\mu'(\vr_n)(\lap\mu(\vr_n))^2}-2\kappa\intO{\mu'(\vr_n)\pt\vr_n\lap\mu(\vr_n)}\\
=\intO{\mu'(\vr_n)(\pt\vr_n)^2}+\kappa^2\intO{\mu'(\vr_n)(\lap\mu(\vr_n))^2}
+\kappa\Dt\intO{|\Grad\mu(\vr_n)|^2}\leq c(t).\label{h_reg}
}
Note that we have passed to the limit with $\alpha$ in the previous paragraph.

Note that this estimate asks for $L^\infty((0,T)\times\Omega)$ bound for $\vw_n$, which is possible only at the level of Galerkin approximation.
 Nevertheless, regularity \eqref{h_reg} allows us to  first rewrite \eqref{A_cont}  as
\eq{\pt\mu(\vr_n) + \Div([\vw_n]_\delta\mu(\vr_n))-2\kappa\mu'(\vr_n)\lap\mu(\vr_n) +\lr{\mu'(\vr_n)\vr_n-\mu(\vr_n)}\Div[\vw_n]_\delta= 0.\label{mu_w}}
Multiplying the above equation by $-\lap\mu(\vr)$ and integrating by parts we obtain
\eqh{&\Dt\intO{|\Grad\mu(\vr_n)|^2}+2\kappa\intO{\mu'(\vr_n)|\lap\mu(\vr_n)|^2}\\
&=\intO{\mu'(\vr_n)\vr_n\Div[\vw_n]_\delta\lap\mu(\vr_n)}+\intO{[\vw_n]_\delta\Grad\mu(\vr_n)\lap\mu(\vr_n)}.}
Note that for any $f\in H^2(\Omega)\cap L^\infty(\Omega)$ we have the following realization of the Gagliardo-Nierenberg inequality
\eq{\|\Grad f\|^2_{L^4(\Omega)}\leq c\|\lap f\|_{L^2(\Omega)}\|f\|_{L^\infty(\Omega)},\label{G-N}}
see for example \cite{GNI}. Therefore, applying this inequality with $f=\mu(\vr_n)$ and using the uniform bound for $\vw_n\in L^2(0,T;H^1(\Omega))$ and $\ep$-dependent bound in $L^4(0,T;L^4(\Omega))$ following from \eqref{ueu},  we obtain
\eq{ \label{est}
\|\Grad\mu(\vr_n)\|_{L^\infty(0,T; L^2(\Omega))}+2\kappa^{1/2}\|\lap\mu(\vr_n)\|_{L^2((0,T)\times\Omega)}\leq c\lr{\|\vr_n^0\|_{H^1(\Omega)},R,\kappa,\ep},
}
and thus, coming back to \eqref{mu_w} and using \eqref{max} we easily check the uniform (with respect to $n$ and $\delta$ ) estimate
\eq{\|\Grad\vr_n\|_{L^\infty(0,T;L^{2}(\Omega))}
+\|\lap\vr_n\|_{L^2((0,T)\times\Omega)}\leq c.
\label{high}}

\bigskip

\noindent{\bf Passage to the limit $n\to\infty$.} As in the previous paragraph the biggest problem is to identify the limit of the  term
\eq{\mu'(\vr_n)\Grad\vr_n\otimes\vw_n,\label{conv_w}}
and in the convective term
\eq{\vr_n\vw_n\otimes\vw_n.\label{conv_w2}}
Having obtained estimate \eqref{high} we can estimate the time-derivative of gradient of $\vr_n$. Indeed taking gradient of  \eqref{A_cont} we obtain
\eqh{\pt\Grad\vr_n=-\Grad\lr{\vr_n\Div[\vw_n]_\delta}-\Grad\lr{[\vw_n]_\delta\cdot\Grad\vr_n}-\kappa\Grad\lap\mu(\vr_n)\in L^2(0,T; W^{-1,3/2}(\Omega)).}
Note, however, that unlike in the analogous step of \cite{BrGiZa} the above estimate is no longer uniform with respect to $\ep$.
Now, applying the Aubin-Lions lemma for $\Grad\vr_n$ we obtain
\eqh{\Grad\vr_n\to\Grad\vr\quad\text{strongly in } L^2(0,T;L^2(\Omega)),}
therefore due to \eqref{max} we also have
\eq{\vr_n\to\vr\quad\text{and}\quad
\frac{1}{\vr_n}\to\frac{1}{\vr}\quad\text{strongly in } L^p(0,T;L^p(\Omega))
\label{s_rho}}
for $p<\infty$ and
\eq{\vr_n\vw_n\to\vr\vw\quad \text{weakly in } L^{p_1}(0,T; L^{q_1}(\Omega))\cap L^{p_2}(0,T; L^{q_2}(\Omega)),\label{k1}}
where $ p_1<2, q_1<6, \ p_2<\infty,q_2<2.$
These convergences justify the limit passage in \eqref{conv_w}.

To justify passage to the limit in \eqref{conv_w2} we first  estimate 
\eq{\|\Grad\lr{\vr_n\vw_n}\|_{L^2(0,T; L^{\frac{3}{2}}(\Omega))}
\leq &\|\Grad\vr_n\|_{L^\infty(0,T;L^2(\Omega))}\|\vw_n\|_{L^2(0,T; L^6(\Omega))}\\
&+\|\Grad\vw_n\|_{L^2(0,T;L^2(\Omega))}\|\vr_n\|_{L^\infty(0,T;L^\infty(\Omega))}.
\label{k2}}
We can also estimate the time derivative of momentum, from \eqref{FG_a} we obtain
\eq{
&\sup_{\|\vcg{\phi}\|\leq1}\left|\intTO{\pt\lr{\vr_n\vw_n}\cdot\vcg{\phi}}\right|\\
&=\sup_{\|\vcg{\phi}\|\leq 1}\left\{\left|\intTO{((\vr_n [\vw_n]_\delta - 2\kappa\Grad\mu({\vr_n}) ) \otimes \vw_n):\Grad \vcg{\phi}}\right|\right.\\
&\quad\qquad+2(1-\kappa)\left|\intTO{\mu(\vr_n)D(\vw_n):\Grad \vcg{\phi}}\right|
+2\kappa\left|\intTO{\mu(\vr_n)A(\vw_n): \Grad\vcg{\phi}}\right|\\
&\quad\qquad+ 2(1-\kappa)\left|\intTO{\lr{\mu'(\vr_n)\vr_n -\mu(\vr_n) }\Div \, \vw_n\, \Div \, \vcg{\phi}} \right|\\
&\quad\qquad+2\kappa(1-\kappa)\left|\intTO{\mu(\vr_n)\Grad\vv_n:\Grad\vcg{\phi}}\right|
+\left|\intTO{p(\vr_n)\Div \, \vcg{\phi}}\right|
\\
&\quad\qquad+2\kappa(1-\kappa)\left| \intTO{(\mu'(\vr_n)\vr_n -\mu(\vr_n) )\Div \, \vv_n\, \Div \, \vcg{\phi}} \right| \\
&\quad\qquad\left.+\ep\left|\intTOB{\lap\vw^s\cdot\lap^s\vcg{\phi}}\right|
+\ep\left|\intTOB{(1+|\Grad\vw|^2) \Grad\vw:\Grad\vcg{\phi}}\right|\right\},
\label{ptw}}
where the norm $\|\vcg{\phi}\|$ denotes the norm in the space $W_T:=L^2(0,T;H^{2s}(\Omega))\cap L^4(0,T;W^{1,4}(\Omega))$ with $s\geq2$.
Let us estimate the convective term
\eqh{&\left|\intTO{((\vr_n [\vw_n]_\delta - 2\kappa\mu'({\vr_n}) \Grad \vr_n) \otimes \vw_n):\Grad \vcg{\phi}}\right|\\
&\leq\intT{\|\Grad\vcg{\phi}\|_{L^6(\Omega)}\lr{R\|\vw_n\|^2_{L^\frac{12}{5}(\Omega)}+c(\kappa,R)\|\Grad\vr_n\|_{L^6(\Omega)}\|\vw_n\|_{L^\frac{3}{2}(\Omega)}}}\\
&\leq c(\kappa,R,\ep)\|\vcg{\phi}\|_{L^2(0,T; W^{2,2}(\Omega))},
}
that is bounded for $s\geq 2$. For the highest order terms we have
\eqh{\ep\left|\intTO{\lap^s\vw_n\cdot\lap^s\vcg{\phi}}\right|\leq \ep\|\vw_n\|_{L^2(0,T; H^{2s}(\Omega))}\|\vcg{\phi}\|_{L^2(0,T; H^{2s}(\Omega)),}}
and
\eqh{\ep\left|\intTO{(1+|\Grad\vw|^2) \Grad\vw:\Grad\vcg{\phi}}\right|
\leq \ep\intT{\|\Grad\vcg{\phi}\|_{L^4(\Omega)}\lr{\|\Grad\vw_n\|^3_{L^{4}(\Omega)}+\|\Grad\vw_n\|_{L^{\frac{4}{3}}(\Omega)}}}\\
\leq \ep\|\Grad\vcg{\phi}\|_{L^4(0,T; L^4(\Omega))}\lr{\|\Grad\vw\|^3_{L^4(0,T; L^4(\Omega))}+\|\Grad\vw\|_{L^{\frac{4}{3}}(0,T; L^{\frac{4}{3}}(\Omega))}}
.}
The other terms from \eqref{ptw} are less restrictive, therefore
\eq{\|\pt\lr{\vr_n\vw_n}\|_{W_T^*}\leq c,\label{k3}}
where $W_T^*$ denotes the dual space of $W_T$ defined above.
Collecting \eqref{k1}, \eqref{k2}, \eqref{k3} and applying the Aubin-Lions lemma to $\vr_n\vw_n$, we prove that
\eqh{\vr_n\vw_n\to\vr\vw\quad \text{strongly in } L^p(0,T;L^p(\Omega))}
for some $p>1$ and therefore thanks to \eqref{s_rho} and \eqref{ueu}
\eq{\Grad\vw_n\to\Grad\vw\quad \text{strongly in }L^p(0,T;L^p(\Omega))\label{strong_gw}}
 for $1\leq p<4$. In particular, convergence in \eqref{conv_w2} is  proved.

For future purposes we now estimate the time derivative of $\vr\vv$. We use \eqref{FGtheta_a}
to obtain
 \eqh{
&\sup_{\|\vcg{\xi}\|\leq1}\left|\intTO{\pt\lr{\vr_n\vv_n}\cdot{\vcg{\xi}}}\right|\\
&= \sup_{\|\vcg{\xi}\|\leq1}\left\{\left|\intTO{((\vr_n [\vw_n]_\delta - 2\kappa \mu'({\vr_n}) \Grad \vr_n) \otimes \vv_n):\Grad \vcg{\xi}}\right|\right.\\
&\qquad+2\kappa\left|\intTO{\mu(\vr_n)\Grad\vv_n:\Grad \vcg{\xi}}\right|
+2\kappa\left|\intTO{(\mu'(\vr_n)\vr_n - \mu(\vr_n))\Div \, \vv_n\ \Div \, \vcg{\xi}}\right|\\
&\qquad\left.+2\left| \intTO{(\mu'(\vr_n)\vr_n - \mu(\vr_n))\Div \, \vw_n\,\Div \, \vcg{\xi}}\right|+2\left|\intTO{\mu(\vr_n)\Grad^t\vw_n:\Grad\vcg{\xi}}\right|\right\}}
 for $\vcg{\xi}\in L^4(0,T;W^{1,4}(\Omega))$. 
 We will only estimate the convective term since it is most restrictive
\eqh{ 
&\left|\intTO{((\vr_n [\vw_n]_\delta - 2\kappa \mu'({\vr_n}) \Grad \vr_n) \otimes \vv_n):\Grad \vcg{\xi}}\right|\\
&\leq\intT{\|\Grad\vcg{\xi}\|_{L^4(\Omega)}\lr{R\|\vw_n\|_{L^4(\Omega)}\|\vv_n\|_{L^2(\Omega)}+c(\kappa,R)\|\Grad\vr_n\|_{L^4(\Omega)}\|\vv_n\|_{L^2(\Omega)}}}\\
&\leq c(\kappa,R)\|\vcg{\xi}\|_{L^4(0,T; W^{1,4}(\Omega))}\|\vv_n\|_{L^2(0,T;L^2(\Omega))}
\lr{\|\vw_n\|_{L^4(0,T; L^4(\Omega))}+\|\lap\vr_n\|^{\frac{1}{2}}_{L^2(0,T; L^2(\Omega))}}
}
thus
\eq{\|\pt\lr{\vr_n\vv_n}\|_{L^{\frac{4}{3}}(0,T; W^{-1,\frac{4}{3}}(\Omega))}\leq c.\label{k4}}
Hence, the limit functions $(\vr,\vw,\vv)=(\vr_\delta,\vw_\delta,\vv_\delta)$ fulfill
\begin{itemize}
\item the continuity equation 
\eq{\pt\vr_\delta+\Div\lr{\vr_\delta[\vw_\delta]_\delta}-2\kappa\lap\mu(\vr_\delta)=0\label{strong1}}
 a.e. in $(0,T)\times\Omega$,
\item the momentum equation
\begin{equation}\label{FGn}
\begin{split}
&\langle\pt\lr{\vr_\delta\vw_\delta},\vcg{\phi}\rangle_{(W_\tau^*, W_\tau)}
-\inttauO{((\vr_\delta [\vw_\delta]_\delta - 2\kappa \Grad \mu(\vr_\delta)) \otimes \vw_\delta):\Grad \vcg{\phi}}\\
&+2(1-\kappa)\inttauO{\mu(\vr_\delta)D(\vw_\delta):\Grad \vcg{\phi}}
+2\kappa\inttauO{\mu(\vr_\delta)A(\vw_\delta): \Grad\vcg{\phi}}\\
&+ 2(1-\kappa)\inttauO{\lr{\mu'(\vr_\delta)\vr_\delta -\mu(\vr_\delta) }\Div \, \vw_\delta\, \Div \, \vcg{\phi}} \\
&-2\kappa(1-\kappa)\inttauO{\mu(\vr_\delta)\Grad\vv_\delta:\Grad\vcg{\phi}}\\
&  -  2\kappa(1-\kappa) \inttauO{(\mu'(\vr_\delta)\vr_\delta -\mu(\vr_\delta) )\Div \, \vv_\delta\, \Div \, \vcg{\phi}}  \\
&-\inttauO{p(\vr_\delta)\Div \, \vcg{\phi}}+\ep\inttauO{\lr{\lap^s\vw_\delta\cdot\lap^s\vcg{\phi}+(1+|\Grad\vw_\delta|^2) \Grad\vw_\delta:\Grad\vcg{\phi}}}=0
\end{split}
\end{equation}
for all $\vcg{\phi}\in W_\tau$ with $\tau\in[0,T]$,
\item the auxiliary equation for $\vv$
\eq{\label{FGthetan}
&
\langle\pt\lr{\vr_\delta\vv_\delta},\vcg{\xi}\rangle_{(L^{\frac{4}{3}}(0,\tau; W^{-1,\frac{4}{3}}(\Omega)), L^4(0,\tau; W^{1,4}(\Omega)))}
\\
&\quad-\inttauO{((\vr_\delta [\vw_\delta]_\delta - 2\kappa  \Grad \mu(\vr_\delta)) \otimes \vv_\delta):\Grad \vcg{\xi}}+2\kappa\inttauO{\mu(\vr_\delta)\Grad\vv_\delta:\Grad \vcg{\xi}}\\
&\quad+2\kappa\inttauO{(\mu'(\vr_\delta)\vr_\delta - \mu(\vr_\delta))\Div \, \vv_\delta\ \Div \, \vcg{\xi}}\\
&-2 \inttauO{(\mu'(\vr_\delta)\vr_\delta - \mu(\vr_\delta))\Div \, \vw_\delta\,\Div \, \vcg{\xi}}-2\inttauO{\mu(\vr_\delta)\Grad^t\vw_\delta:\Grad\vcg{\xi}}=0,
}
for $\vcg{\xi}\in L^4(0,\tau; W^{1,4}(\Omega))$, $\tau\in[0,T]$.
\end{itemize}

\subsection{Passage to the limit $\delta$ tends to $0$ and identification of $\vv$ with $2\Grad\vp(\vr)$}
The aim of this paragraph is to let $\delta\to0$ in the equations \eqref{strong1}, \eqref{FGn} and  \eqref{FGn}. This limit passage can be performed exactly as $n\to \infty$ presented above. The only difference is that after this step we may drop the additional equation for $\vv$ thanks to identification $\vv=2\Grad\vp(\vr)$. Below we present the details of this reasoning.

Note that the coefficients of the quasi-linear parabolic equation \eqref{strong1} (i.e. $[\vw_\delta]_\delta$) are sufficiently regular and the maximum principle \eqref{max} holds uniformly with respect to all approximation parameters except $\varepsilon$. Therefore, the classical theory  of Lady{\v{z}}enskaja, Solonnikov and Uralceva \cite{LSU} (Theorems 7.2, 7.3 and 7.4 from \cite{LSU}) can be applied to show further regularity of $\vr_\delta$, we have in particular
\eq{\pt\vr_\delta\in C([0,T]; C(\Omega)),\quad \vr_\delta\in C([0,T]; C^2(\Omega)).\label{reg_rg}}
Let us now rewrite \eqref{strong1} using \eqref{vp_mu} as
$$  \pt\vr_\delta
     + \Div\lr{\vr_\delta  ([\vw_\delta]_\delta - 2\kappa\, \Grad\vp(\vr_\delta)) } = 0$$  
and therefore multiplying the above equation by $\mu'(\vr)$ we obtain 
$$\partial_t\mu(\vr_\delta) 
     + \Div\lr{\mu(\vr_\delta)([\vw_\delta]_\delta 
       - 2\kappa\, \Grad\vp(\vr_\delta))} 
     +(\mu'(\vr_\delta)\vr_\delta - \mu(\vr_\delta))\lr{\Div[\vw_\delta]_\delta-2\kappa \Delta \vp(\vr_\delta) } 
     = 0. $$  
Differentiating it with respect to space one gets in the sense of distributions
\eq{ \pt(\vr_\delta \widetilde \vv_\delta) 
     + \Div(\vr_\delta ([\vw_\delta]_\delta 
     - 2\kappa\, \Grad\vp(\vr_\delta))\otimes \widetilde\vv_\delta ) 
    +2\Grad\lr{(\mu'(\vr_\delta)\vr_\delta - \mu(\vr_\delta))\lr{\Div[\vw_\delta]_\delta-\kappa \Div \, \tilde\vv_\delta}}\\
    + 2\Div((\mu(\vr_\delta) \Grad^{t} [\vw_\delta]_\delta)
    - 2\kappa \Div(\mu(\vr_\delta)\Grad\widetilde \vv_\delta)
      = 0\label{test}}
 where by $\widetilde\vv_\delta$ we denoted $2\Grad\vp(\vr_\delta)$. 
 Note that due to particular case of the Gagliardo-Nirenberg interpolation inequality \eqref{G-N} and \eqref{max} we know that $\Grad\vr_\delta$ is bounded
  \eqh{\|\Grad \vr_\delta\|_{L^4(0,T; L^4(\Omega))}\leq c,} 
 uniformly with respect to $\delta$.  One can thus estimate the convective term of \eqref{test} in $L^2(0,T; W^{-1,2}(\Omega))$ uniformly with respect to $\delta$. Indeed, we now $\vw_\delta$ uniformly
bounded in $L^4(0,T;L^4(\Omega))$ with respect to $\delta$ and therefore
\eqh{&\sup_{\|\vcg{\xi}\|\leq1}\left|\intTO{  (\vr_\delta ([\vw_\delta]_\delta 
     - 2\kappa\, \Grad\vp(\vr_\delta))\otimes \widetilde\vv_\delta ):\Grad\vcg{\xi}}\right|\\
     &\leq c(R) \|\Grad\vcg{\xi}\|_{L^2(0,T; L^2(\Omega))}\|\Grad \vr_\delta\|_{L^4(0,T; L^4(\Omega))}
     \lr{\|\Grad [\vw_\delta]_\delta\|_{L^4(0,T; L^4(\Omega))}+\|\Grad \vr_\delta\|_{L^4(0,T; L^4(\Omega))}}
     }
 for $\vcg{\xi}\in L^2(0,T; W^{1,2}(\Omega))$ (uniformly with respect to $\delta$), which justifies that 
\eq{ &\langle \pt(\vr_\delta \widetilde \vv_\delta), \vcg{\xi}\rangle_{(L^2(0,T; W^{-1,2}(\Omega)),L^2(0,T; W^{1,2}(\Omega)))}
-\intTO{\vr_\delta ([\vw_\delta]_\delta 
     - 2\kappa\, \Grad\vp(\vr_\delta))\otimes \widetilde\vv_\delta :\Grad\vcg{\xi}}\\
     &-2\intTO{ (\mu'(\vr_\delta)\vr_\delta - \mu(\vr_\delta))  \Div [\vw_\delta]_\delta\, \Div \, \vcg{\xi} }+2\kappa\intTO{ (\mu'(\vr_\delta)\vr_\delta - \mu(\vr_\delta))  \Div \, \widetilde\vv_\delta\, \Div \, \vcg{\xi} }\\
   & -2\intTO{\mu(\vr_\delta) \Grad^{t} [\vw_\delta]_\delta:\Grad\vcg{\xi}}
    +2\kappa \intTO{\mu(\vr_\delta)\Grad\widetilde \vv_\delta:\Grad\vcg{\xi}}
      = 0\label{test2}}
 is satisfied for any $\vcg{\xi}\in L^2(0,T; W^{1,2}(\Omega))$.
 \begin{rmk}
 As we noticed in \eqref{reg_rg} the regularity of $\vr_\delta$ is in fact much higher and allows to formulate the equation for $\widetilde\vv_\delta$ in much stronger sense than merely \eqref{test2}. This formulation, however, will be used when passing to the limit with respect to $\ep$ after passing to the limit with respect to $\delta$.
 \end{rmk}
 We now want  show that  $\widetilde\vv_\delta -\vv_\delta$ tends to $0$ when $\delta$ goes to zero
 in an appropriate norm. To this purpose let us expand
  \eq{
I= & \Dt\intO{\vr_\delta\frac{|\vv_\delta-\widetilde\vv_\delta|^2}{2}}
  + 2\kappa \intO{\mu(\vr_\delta)|\Grad (\vv_\delta - \widetilde\vv_\delta)|^2} \\
  & +2\kappa\intTO{(\mu'(\vr_\delta)\vr_\delta - \mu(\vr_\delta))(\Div(\vv_\delta -\widetilde \vv_\delta))^2}\\
 =&  \Dt\intO{\vr\lr{\frac{|\vv_\delta|^2}{2}-\vv_\delta \cdot\widetilde\vv_\delta+\frac{|\widetilde\vv_\delta|^2}{2}}}\\ 
 &+2 \kappa\intO{\mu(\vr_\delta)\lr{|\Grad\vv_\delta |^2 +|\Grad \widetilde\vv_\delta|^2
    - 2 \Grad \vv_\delta\cdot \Grad \widetilde \vv_\delta}}. \\
 & +2\kappa\intTO{(\mu'(\vr_\delta)\vr_\delta - \mu(\vr_\delta))\lr{|\Div \, \vv_\delta|^2 + |\Div \, \widetilde \vv_\delta|^2  - 2 \Div \, \vv_\delta \Div \, \widetilde \vv_\delta}}.
  \label{expand}
 }
 To handle the first term, let us notice that letting $n\to \infty$ in \eqref{t1}, using the lower semi-continuity of the convex functions and the strong convergence of $\Grad\vw_n$ established in \eqref{strong_gw} we obtain
 \eq{\label{t11}
&\Dt\intO{\vr_\delta\frac{|\vv_\delta|^2}{2}}
+2\kappa\intTO{\mu(\vr_\delta)|\Grad\vv_\delta|^2}+2\kappa\intTO{(\mu'(\vr_\delta)\vr_\delta - \mu(\vr_\delta))(\Div \, \vv_\delta)^2}\\
&\quad-2 \intO{(\mu'(\vr_\delta)\vr_\delta - \mu(\vr_\delta))\Div \, \vw_\delta\,\Div \, \vv_\delta}-2\intTO{\mu(\vr_\delta)\Grad^t\vw_\delta:\Grad\vv_\delta}\leq0.
}
 Now, the last term in \eqref{expand} can be computed using $\vcg{\xi}=\widetilde\vv_\delta$ in \eqref{test2}, we have
  \eqh{&\Dt\intO{\vr_\delta\frac{|\widetilde\vv_\delta|^2}{2}} 
  +2\kappa \intTO{\mu(\vr_\delta)|\Grad\widetilde \vv_\delta|^2}
  +2\kappa\intO{ (\mu'(\vr_\delta)\vr_\delta - \mu(\vr_\delta))  (\Div \, \widetilde\vv_\delta)^2 }\\
   & -2\intO{\mu(\vr_\delta) \Grad^{t} [\vw_\delta]_\delta :\Grad\widetilde\vv_\delta}
    -2\intTO{ (\mu'(\vr_\delta)\vr_\delta - \mu(\vr_\delta))  \Div \, \vw\, \Div \, \tilde\vv_\delta }=0.
 }
The middle term in \eqref{expand} equals
\eq{\Dt\intO{\vr_\delta\vv_\delta\cdot\widetilde\vv_\delta}=\intOB{\pt\lr{\vr_\delta\vv_\delta}\cdot\widetilde\vv_\delta + \vv_\delta \cdot\pt(\vr_\delta \widetilde\vv_\delta ) - \partial_t\vr_\delta\,  \vv_\delta\cdot \widetilde \vv_\delta}\label{der_mixed}}
and the two first terms make sense and can be handled using $\vcg{\xi}=\widetilde\vv_\delta$ in \eqref{FGthetan} and $\vcg{\xi} = \vv_\delta$ in \eqref{test2}. 
Note that $\widetilde \vv_\delta$ and $\partial_t \vr_\delta$ are due to \eqref{reg_rg} regular enough to justify the integrability of the last term in \eqref{der_mixed} and we can write
\eqh{&\intO{\pt\vr_\delta\vv_\delta\cdot\widetilde\vv_\delta}\\
&=\intO{\lr{\vr_\delta  ([\vw_\delta]_\delta - 2\kappa\, \Grad\vp(\vr_\delta)) } \otimes\vv_\delta:\Grad\widetilde\vv_\delta}
+\intO{\lr{\vr_\delta  ([\vw_\delta]_\delta - 2\kappa\, \Grad\vp(\vr_\delta)) } \otimes\widetilde\vv_\delta:\Grad\vv_\delta}.
} 
Therefore, after summing all expressions together and after some manipulation, we can show that
\eqh{I
&-2\intTO{ (\mu'(\vr_\delta)\vr_\delta - \mu(\vr_\delta))  \Div [\vw_\delta]_\delta\, \Div \, \tilde\vv_\delta}
+2\intTO{ (\mu'(\vr_\delta)\vr_\delta - \mu(\vr_\delta))  \Div \, \vw_\delta\, \Div \, \tilde\vv_\delta}\\
&+2\intTO{ (\mu'(\vr_\delta)\vr_\delta - \mu(\vr_\delta))  \Div [\vw_\delta]_\delta\, \Div \, \vv_\delta }
-2\intTO{ (\mu'(\vr_\delta)\vr_\delta - \mu(\vr_\delta))  \Div \, \vw_\delta\, \Div \, \vv_\delta }\\
&-\intO{\mu(\vr_\delta)\Grad^t[\vw_\delta]_\delta:\Grad\widetilde\vv_\delta}
-\intO{\mu(\vr_\delta)\Grad^t\vw_\delta:\Grad\vv_\delta}\\
&+\intO{\mu(\vr_\delta)\Grad^t\vw_\delta:\Grad\widetilde\vv_\delta}
+\intO{\mu(\vr_\delta)\Grad^t[\vw_\delta]_\delta:\Grad\vv_\delta}\leq 0,
}
in particular
\eqh{I \leq &\intO{\mu(\vr)\Grad^t([\vw_\delta]_\delta - \vw_\delta ): \Grad (\widetilde\vv_\delta - \vv_\delta)}\\
&+
2\intTO{ (\mu'(\vr_\delta)\vr_\delta - \mu(\vr_\delta))  \Div \lr{[\vw_\delta]_\delta-\vw_\delta} \Div\lr{\tilde\vv_\delta-\vv_\delta}}.}
Note that the r.h.s. of this inequality tends to 0 when $\delta\to 0$. Indeed, we can bound $ \Grad (\widetilde\vv_\delta - \vv_\delta)$ in $L^2(0,T; L^2(\Omega))$ uniformly with respect to $\delta$ and $[\vw_\delta]_\delta\to \vw$ strongly in $L^p(0,T; W^{1,p}(\Omega))$ for $p<4$.
Therefore, using \eqref{expand}, we conclude that $\vv_\delta - \widetilde\vv_\delta$ converges to zero in $L^\infty(0,T;L^2(\Omega))\cap L^2(0,T;H^1(\Omega))$ when $\delta \to 0$. $\Box$

The limit functions $(\vr,\vw)=(\vr_\ep,\vw_\ep)$ fulfill
\begin{itemize}
\item the continuity equation 
\eq{\pt\vr_\ep+\Div\lr{\vr_\ep\vw_\ep}-\kappa\lap\mu(\vr_\ep)=0\label{cont_ep}}
 a.e. in $(0,T)\times\Omega$,
\item the momentum equation
\begin{equation}\label{FG_ep}
\begin{split}
&%-\intO{\vr\vw\cdot\pt{\bf{\vcg{\phi}}}}
\langle\pt\lr{\vr_\ep\vw_\ep},\vcg{\phi}\rangle_{(W_\tau^*, W_\tau)}
-\inttauO{((\vr_\ep\vw_\ep - 2\kappa \Grad \mu(\vr_\ep)) \otimes \vw_\ep):\Grad \vcg{\phi}}\\
&\quad+2(1-\kappa)\inttauO{\mu(\vr_\ep)D(\vw_\delta):\Grad \vcg{\phi}}+2\kappa\intTO{\mu(\vr_\ep)A(\vw_\ep): \Grad\vcg{\phi}}\\
&\quad+ 2(1-\kappa)\inttauO{\lr{\mu'(\vr_\ep)\vr_\ep -\mu(\vr_\ep) }\Div \, \vw_\ep\, \Div \, \vcg{\phi}} \\
&\quad-4\kappa(1-\kappa)\inttauO{\mu(\vr_\ep)\Grad^2\vp(\vr_\ep):\Grad\vcg{\phi}}\\
&  \quad-  4\kappa(1-\kappa) \inttauO{(\mu'(\vr_\ep)\vr_\ep -\mu(\vr_\ep) )\lap\vp(\vr_\ep)\, \Div \, \vcg{\phi}}  \\
&\quad-\inttauO{p(\vr_\ep)\Div \, \vcg{\phi}}+\ep\inttauO{\lr{\lap^s\vw_\ep\cdot\lap^s\vcg{\phi}+(1+|\Grad\vw_\ep|^2) \Grad\vw_\ep:\Grad\vcg{\phi}}}=0
\end{split}
\end{equation}
for $\vcg{\phi}\in W_\tau$,  $\tau\in[0,T]$,
\item the auxiliary equation for $\Grad\vp(\vr_\ep)$
\eq{\label{FGtheta_ep}
&
\langle\pt\Grad\mu(\vr_\ep),\vcg{\xi}\rangle_{(L^2(0,\tau; W^{-1,2}(\Omega)),L^2(0,\tau; W^{1,2}(\Omega)))}
\\
&\quad-\inttauO{((\vr_\ep\vw_\ep- 2\kappa  \Grad \mu(\vr_\ep)) \otimes \Grad\vp(\vr_\ep)):\Grad \vcg{\xi}}+2\kappa\inttauO{\Grad\mu(\vr_\ep):\Grad \vcg{\xi}}\\
&\quad+2\kappa\inttauO{(\mu'(\vr_\ep)\vr_\ep - \mu(\vr_\ep))\lap\vp(\vr_\ep)\ \Div \, \vcg{\xi}}\\
&\quad- \inttauO{(\mu'(\vr_\ep)\vr_\ep - \mu(\vr_\ep))\Div \, \vw_\ep\,\Div \, \vcg{\xi}}
-\inttauO{\mu(\vr_\ep)\Grad^t\vw_\ep:\Grad\vcg{\xi}}=0,
}
for all  $\vcg{\xi}\in L^2(0,\tau; W^{1,2}(\Omega))$ with $\tau\in[0,T]$.
\end{itemize}

\subsection{Recovering the global weak $\kappa$-entropy solution}\label{SS:Recov}
 Due to the identification, we may now define a new quantity
\eq{\vu=\vw- 2\kappa\Grad\vp(\vr).}
 Compared to \cite{BrGiZa}, we cannot pass to the limit in the equations satisfied by $\vw$ and $\vv$. We 
 have to combine them to deal with an equation on $\vu$. More precisely, multiplying \eqref{FGtheta} by $\kappa$ and substracting to \eqref{FG} we obtain the weak formulation on $\vu$
\begin{equation}\label{u}
\begin{split}
&{\langle\pt(\vr_\ep\vu_\ep),\vcg{\phi}\rangle}_{(W_\tau^*,W_\tau)}
-\inttauO{(\vr_\ep \vu_\ep \otimes \vu_\ep):\Grad \vcg{\phi}}
+2\inttauO{\mu(\vr_\ep)D(\vu_\ep):\Grad \vcg{\phi}}\\
&\quad +\inttauO{\lambda(\vr_\ep)\Div \, \vu_\ep\ \Div \, \vcg{\phi}}
-\inttauO{p(\vr_\ep)\Div \, \vcg{\phi}}\\
&\quad+\ep\inttauO{\lap^{s}\vw_\ep\cdot\lap^{s}\vcg{\phi}} + \ep\inttauO{(1+|\Grad\vw_\ep|^2)\Grad \vw_\ep \cdot\Grad\vcg{\phi}} =0 
\end{split}
\end{equation}
satisfied for all $\vcg{\phi}\in W_\tau$ with $\tau\in[0,T]$.\\
Passing to the limit in \eqref{u} with respect to $\varepsilon$ follows the lines introduced by \cite{BrDe4}
using the uniform estimates given by the regularized $\kappa$-entropy.
See also \cite{Za} for details regarding the barotropic case with chemical species.
Now, after passing to the limit in the energy estimate \eqref{aa1} (with all parameters $\alpha\to0$, $n\to\infty$, $\delta\to0$ in the same way), we obtain
 \eq{\label{aaa}
&\Dt\intO{\vr_\ep\lr{\frac{|\vw_\ep|^2}{2}
+  ( 1-\kappa)\kappa\frac{|2 \Grad\vp(\vr_\ep)|^2}{2}  }}
 +  \Dt \intO{ \vr_\ep e(\vr_\ep)} + 2(1-\kappa)\intO{\mu(\vr_\ep) |D(\vu_\ep)|^2} \\
&\quad+ 2 \kappa\intO{\mu(\vr_\ep) |A(\vw_\ep)|^2}+ 2 (1-\kappa)\intO{(\mu'(\vr_\ep)\vr_\ep-\mu(\vr_\ep))(\Div \, \vu_\ep)^2}
\\
&\quad+ 2\kappa \intO {\frac{\mu'(\vr_\ep)p'(\vr_\ep)}{\vr} |\Grad\vr_\ep|^2} 
+ \ep\intO {|\Delta^{s} \vw_\ep|^2}
+ \ep \intO {(1+|\Grad\vw_\ep|^2)|\Grad\vw_\ep|^2} \leq 0,
}
where we essentially used the regularization  $\ep [\lap^{2s} \vw_\ep - \Div((1+|\Grad \vw_\ep|^2)\Grad \vw_\ep)$.
We then  show that the limit when $\ep\to$ give a solution $(\vr,\vu)$ which satisfies the following energy inequality
\eq{\label{energy}
&\sup_{t\in [0,T] }\Bigl[\intO{\vr\lr{\frac{|\vw|^2}{2}
+ ( 1-\kappa)\kappa\frac{|2 \Grad\vp(\vr)|^2}{2}  }(t)}
 +  \intO{ \vr e(\vr)(t)}\Bigr] +2(1-\kappa)\int_0^T\intO{\mu(\vr) |D(\vu)|^2} \\
&+2 \kappa\int_0^T\intO{\mu(\vr) |A(\vw)|^2}
+2 (1-\kappa)\int_0^T\intO{(\mu'(\vr)\vr-\mu(\vr))(\Div \, \vu)^2}
\\
&+ 2 \kappa \int_0^T\intO {\frac{\mu'(\vr)p'(\vr)}{\vr} |\Grad\vr|^2} 
 \le \intOB{\vr_0\lr{\frac{|\vw_0|^2}{2}
+ ( 1-\kappa)\kappa\frac{|2 \Grad\vp(\vr_0)|^2}{2}  }
 +  \vr_0 e(\vr_0)} .
}
For the limit passage in the corresponding equations we use the density estimates due to the singular pressure, as it was introduced in the work by the  first two authors and generalized to
chemically reacting mixture case by the third author and collaborators.
 In particular, one needs to ensure that
$$\vr^{\frac{1}{2}}\vu\in L^s(0,T; L^r(\Omega))$$
with $r,s>2$. Here assumptions of $n,m$ and $\gamma^-$ \eqref{mn}, \eqref{gpm}
play an important role. Indeed, assuming for the moment that $\vr^{\frac{1}{2}}\vu\in L^\infty(0,T; L^2(\Omega))$, $\vr\in L^\infty(0,T; L^p(\Omega))$ and $\vu\in L^{q_1}(0,T; L^{q_2}(\Omega))$, 
we can write exactly as in the proof of Lemma 6.1. from \cite{BrDe3} that
$$\vr^{\frac{1}{2}}|\vu|=\vr^{\frac{1}{2}-\vt}\vr^\vt |\vu|^{2\vt}|\vu|^{1-2\vt}$$
for some $0\leq\vt\leq\frac{1}{2}$. Therefore
$$\|\vr^{\frac{1}{2}}\vu\|_{L^s(0,T; L^r(\Omega))}\leq\|\vr\|^{\frac{1}{2}-\vt}_{L^\infty(0,T; L^p(\Omega))}
\|\vr^{\frac{1}{2}}\vu\|^{2\vt}_{L^\infty(0,T; L^2(\Omega))}\|\vu\|^{1-2\vt}_{L^{q_1}0,T; L^{q_2}(\Omega)}$$
with
\eq{\frac{1}{s}=\frac{1-2\vt}{q_1}, \quad \frac{1}{r}=\frac{1}{p}\lr{\frac{1}{2}-\vt}+\vt+\frac{1-2\vt}{q_2}.\label{sr}}
It follows from the $\kappa$-entropy estimate \eqref{energy} that the best $p$ we can take is equal
$$p=6m-3.$$
Let us now determine $q_1, q_2$,  we first check that the $\kappa$-entropy estimate \eqref{energy} gives us uniform bound on 
$$\intTO {\frac{\mu'(\vr)p'(\vr)}{\vr} |\Grad\vr|^2} .$$
In particular, for $\vr<\vr^*$, we obtain
$$\Grad\lr{\vr^{-\frac{\gamma^-+1-n}{2}}}\in L^2(0,T;L^2(\Omega)),$$
which controls the negative powers of $\vr$ close to vacuum. We can use this estimate to determine $q_1$ and some $q_3$, we have
\eq{\|\Grad\vu\|_{L^{q_1}(0,T; L^{q_3}(\Omega))}\leq c\lr{1+\|\vr^{-\frac{n}{2}}\|_{L^{2j}(0,T; L^{6j}(\Omega))}}
\|\vr{\frac{n}{2}}\Grad\vu\|_{L^2(0,T; L^2(\Omega))}\\
\leq c\lr{1+\|\Grad\lr{\vr^{-\frac{\gamma^-+1-n}{2}}}\|_{L^{2}(0,T; L^{6}(\Omega))}}
\|\vr{\frac{n}{2}}\Grad\vu\|_{L^2(0,T; L^2(\Omega))}
}
where
\eq{j=\frac{\gamma^-+1-n}{n},\quad q_1=\frac{2j}{j+1}=2\lr{1-\frac{n}{\gamma^-+1}}, \quad
\frac{1}{q_3}=\frac{1}{6j}+\frac{1}{2}. 
\label{q1q2}}
Thus in \eqref{sr} we can take 
$$q_2=\frac{3q_3}{3-q_3}=3 q_1$$
provided that
\eqh{q_1>1,\quad \text{i. e.}\quad n<\frac{\gamma^-+1}{2}.}
Inserting $q_2=3q_1$ and $p=6m-3$ to \eqref{sr} and using \eqref{q1q2} we obtain $r,s>2$ provided that in addition we assume that
$$\gamma^->\frac{2n(3m-2)}{4m-3}-1, \quad m>\frac{3}{4},$$
but there is no restriction on $\gamma^+>1$.

\section{Navier--Stokes equations with drag terms} \label{drag}

Let us recall the compressible Navier Stokes system with quadratic turbulent drag force
      \begin{equation}\label{main3}
     \left\{ 
      \begin{array}{l}
            \vspace{0.2cm}
       \pt\vr+\Div\lr{\vr\vu}= 0,\\
            \vspace{0.2cm}
        \pt\lr{\vr\vu} + \Div ({\vr\vu\otimes\vu})  - \Div\lr{2 \mu(\vr) D(\vu)} 
           - \Grad (\lambda(\vr) \Div \, \vu) +  r_1\vr |u|u +  \Grad p(\vr) = 0,\\
                      \end{array}\right. 
    \end{equation}
 with the usual pressure law $p(\vr)= a\vr^\gamma$ with $\gamma>1$.
The augmented regularized system reads
\eq{\label{approx3}
&\pt\vr+ \Div(\vr [\vw]_\delta) - 2 \kappa\Div\lr{[\mu'( \vr)]_\alpha \Grad\vr} = 0,\\
 &\pt\lr{\vr \vw}+ \Div  ((\vr [\vw]_\delta - 2 \kappa  [\mu'( \vr)]_\alpha \Grad\vr) \otimes \vw)  -  \Grad (\lr{\lambda(\vr)- 2 \kappa (\mu'(\vr)\vr-\mu(\vr))}\Div(\vw- \kappa\vv)) \\
&\qquad\qquad- 2(1-\kappa)  \Div (\mu(\vr) D(\vw)) - 2 \kappa \Div (\mu(\vr) A(\vw))  \\
& + \varepsilon \lap^{2s} \vw - \varepsilon \Div ((1+|\Grad\vw|^2)\Grad\vw)  
 + r_1\vr  |\vw-2 \kappa\nabla\varphi(\vr)|(\vw-2 \kappa\varphi(\vr)) + \Grad p(\vr)\\
&\qquad\quad = -  2 \kappa (1-\kappa) \Div(\mu(\vr) \Grad\vv),\\
&\pt(\vr \vv) + \Div((\vr [\vw]_\delta - 2 \kappa [\mu'( \vr)]_\alpha \Grad\vr)\otimes \vv) \\
&\qquad\qquad -2 \kappa \Div(\mu(\vr)\Grad\vv) + 2\Grad \lr{(\mu'(\vr)\vr - \mu(\vr))\Div \lr{\vw-\kappa\vv}}\\
&\qquad\quad= - 2\Div(\mu(\vr) \Grad^{t} \vw),
}
The additional drag term $r_1 |\vw-2\kappa\Grad\varphi(\vr)|(\vw-2\kappa\Grad\varphi(\vr))$  will not really affect the $\kappa$-entropy. More precisely the contribution to the $\kappa$--entropy estimate
may be written as
$$ \intO{\vr |\vw-\kappa\nabla \varphi(\vr)|(\vw-\kappa \nabla \varphi(\vr)) \cdot \vw} = 
  \intO{\vr|\vw-\kappa\nabla \varphi(\vr)|^3} 
    +\kappa \intO{\vr |\vw-\kappa\nabla \varphi(\vr)|(\vw-\kappa \nabla \varphi(\vr)) \cdot \nabla\varphi(\vr)}  $$
    Note that the second term of the above right hand side rewrites:
 $$\int_\Omega \vr |\vw-2\kappa\nabla\varphi(\vr)|(\vw-2\kappa\nabla\varphi(\vr) )\cdot \nabla\varphi(\vr) $$
 $$= - \intO {\mu(\vr) |\vw-2\kappa\nabla\varphi(\vr)| \Div(\vw-2\kappa\nabla\varphi(\vr))}
 - \intO{\mu(\vr) \frac{\vw-2\kappa\nabla\varphi(\vr)}{|\vw-2\kappa\nabla\varphi(\vr)|}\cdot (
  \vw-2\kappa\nabla\varphi(\vr) )\cdot \nabla (\vw-2\kappa\nabla\varphi(\vr))}
 $$
 $$= - \intO {\mu(\vr) |\vu| \Div \, \vu} - \intO{\mu(\vr) \frac{u_k}{|\vu|}  u_j \partial_j u_k}
 $$
where $\vu$ is defined by  $\vu=\vw-2\kappa\Grad\varphi(\vr)$. One may write the estimate 
 $$ \intO{\mu(\vr)|\vu| | D( \vu ) |} \le 
    \|\sqrt{\mu(\vr)}  D(\vu) \|_{L^2(\Omega)} \|\sqrt{\mu(\vr)} \vu\|_{L^2(\Omega)}
    \le  \|\sqrt{\mu(\vr)} D(\vu) \|_{L^2(\Omega)} \|\frac{\sqrt{\mu(\vr)}}{\vr^{1/3}} \vr^{1/3} \vu\|_{L^2(\Omega)}
 $$
 Now we use the hypothesis made on viscosities.  We can bound 
 $$2 r_1\kappa \intO{\mu(\vr)|\vu| |D(\vu)|} \le 
    \frac{1}{2} \|\sqrt{\mu(\vr)} D(\vu)\|^2_{L^2(\Omega)} 
    + \frac{r_1}{3} \|\vr  |\vu|^3\|_{L^1(\Omega)}
    + c(r_1) \intO{\frac{\mu(\vr)^3}{\vr^2}}.
$$
and using \eqref{rhopetit} and \eqref{assum3}
$$ \intO{\frac{\mu(\vr)^3}{\vr^2}} \le C + \intO{\frac{\mu(\vr)^3}{\vr^2}1_{\{(t,x): \vr(t,x)\ge 1 \}}}
 \le C +  c(r_1) \intO{\vr e(\vr)}$$
 with $c(r_1) \to 0$ when $r_1$ tends to zero.
Thus the $\kappa$-entropy is not perturbed and we get uniform estimate on 
 $\vr |\vu|^3$ with respect to $\varepsilon$ in $L^1(0,T;L^1(\Omega))$. 
 It is enough to conclude as in \cite{BrDeGe} because the extra estimate we
 get replaces the Mellet-Vasseur quantity involved in \cite{MeVa}.

\section{Hypocoercivity revisited on linearized compressible Navier-Stokes} \label{hypoc}

   In this section, we want to show to readers who are not familiar with compressible
Navier-Stokes equations why our $\kappa$-entropy equality may be seen as a nonlinear
version of the so-called {\it  hypocoercivity property} which is strongly  used in the  framework of
strong solutions. The interested readers is referred to  \cite{Vi} and \cite{BeZu} for more 
general discussions around hypocoercivity and to  \cite{Da} for  deep mathematical
results on compressible Navier-Stokes equations in critical spaces.

\subsection{Linearized barotropic system}  In this subsection we consider the barotropic
compressible Navier-Stokes system. Linearization around the state $(\vr^0,u^0)=(1,0)$ on the barotropic system gives
\begin{equation}\label{CNSbaro}
\left\{
\begin{array}{l}
\partial_t\vr  + \Div \, \vu =0,\\
\displaystyle \partial_t  \vu + \Grad \vr  - 2\mu \Div  (D(\vu)) +\frac{2}{d} \mu \Grad \Div \, \vu  = \vc{0},\\
\end{array}
\right.
\end{equation}
with $D(\vu) = (\Grad \vu + \Grad^t \vu)/2$.
The chosen coefficient $\lambda$ is the borderline case that means the one such that
$\lambda + {2\mu }/{d} = 0$.
  We focus on the periodic domain  $\Omega=\mathbb{T}^d$.
 The standard energy (multiplying the first equation of \eqref{CNSbaro} by $\vr$ and the second equation by $\vu$ and summing them up) reads
\begin{equation} \label{energybaro}
 \frac{1}{2}\Dt \intOB{ |\vu|^2 + |\vr|^2}
     + \mu \intO{ \left|D(\vu) - \frac{\Div \, \vu}{d} \vc{I}\right|^2}=0.
\end{equation}
Let $\kappa$ be a constant that will be determined later one, then $\vv= \vu+2\kappa\mu \Grad\vr$ satisfies the equation
$$\partial_t \vv + \Grad \vr  - \mu \lap \vu - \mu \lr{1-\frac{2}{d}-2\kappa} \Grad \Div \, \vu = \vc{0}$$
and therefore (multiplying this equation by $\vv$, the first mass equation by $\vr$ and summing 
them up)
\eqh{\frac{1}{2}\Dt \intOB{ |\vv|^2+|\vr|^2 }&+ 2\kappa\mu  \int |\Grad\vr|^2
      + \mu \intO{ |\Grad \vu|^2} + \mu \lr{1- \frac{2}{d} -2\kappa} \intO{ |\Div \, \vu|^2 } \\
  &- 4 \mu^2\kappa \lr{1-\frac{1}{d} - \kappa} \intO{ \Grad \Div \, \vu \cdot \Grad \vr} = 0.}
This equality may be rewritten (using that $\Delta = \Grad \Div - {\rm curl}\, {\rm curl}$) as 
\eq{\label{entropybaro}
&\frac{1}{2}\Dt \intOB{ |\vv|^2+|\vr|^2 }+ 2 \kappa\mu  \intO{ |\Grad\vr|^2}
      +2 \mu \intO{ |A(\vu)|^2} \\
      &+ 2 \mu \lr{1- \frac{1}{d} -\kappa} \intO{ |\Div \, \vu|^2}
 - 4 \mu^2\kappa \lr{1-\frac{1}{d} - \kappa} \intO{ \Grad \Div \, \vu \cdot \Grad \vr} = 0}
 where $A(\vu)=(\nabla\vu -(\nabla u)^t)/2$.
Testing the gradient of the mass equation by $\Grad\vr$, we get
\begin{equation*}
\frac{1}{2} \Dt \intO{ |\Grad\vr|^2} + \intO{ \Grad \vr \cdot \Grad \Div \, \vu} = 0.
\end{equation*}
 Assume now that $0<\kappa < (d-1)/d$, then multiplying the last 
relation by $4 \mu^2\kappa \lr{1-\frac{1}{d} - \kappa}$ and adding to 
\eqref{entropybaro} we get
\eq{\label{entropybarofinal} \displaystyle
&\frac{1}{2}\Dt \intOB{ |\vv|^2+|\vr|^2 }  +2  \mu^2\kappa \displaystyle \lr{1-\frac{1}{d} - \kappa} \Dt \intO{ |\Grad\vr|^2 }\\
&+ 2 \kappa\mu  \intO{ |\Grad\vr|^2}
      +  2\mu \intO{ |A(\vu)|^2} 
      + 2 \mu \lr{1- \frac{1}{d} -\kappa} \intO{ |\Div \, \vu|^2} 
= 0.}
Remark that this equality is the analog, deleting the term $1/d$, of the $\kappa$-entropy we found in the nonlinear framework: Note that
$$\intO{|A(\vu)|^2} + (1-\kappa)\intO{|\Div \, \vu|^2} = \kappa \intO{[A(\vu)|^2}+ (1-\kappa)\intO{|D(\vu)|^2}
$$
 In particular, this $\kappa$-entropy \eqref{entropybarofinal} provides the exponential decay
in time  to $(0,0)$ of $(\vr,\vu)$ in the $L^2(\Omega)$ norm if the initial perturbation
$(\vr_0,\vu_0)$ is in $L^2(\Omega)$.
 
 \bigskip
 Let us now look at exponential decay in$H^1$ norm.
 Taking the $\Curl$ of the momentum equation in \eqref{CNSbaro} and testing against $\Curl  \vu$, we get
 \begin{equation}\label{curlbaro}
 \frac{1}{2}\Dt\intO{ |\Curl  \vu|^2 }+ \mu \intO{ |\Grad \Curl  \vu|^2}= 0.
 \end{equation}
Taking the $\Div$ of the momentum equation and testing it against $\Div \, \vu$
$$ \frac{1}{2}\Dt \intO{ |\Div \, \vu|^2 }+  (2\mu+\lambda) \intO{ |\Grad \Div \, \vu|^2}
       + \intO{ \lap \vr\, \Div \, \vu} = 0.$$
Adding these equation to the mass equation tested against $\Delta\vr$,
we obtain 
\begin{equation}\label{divgradrhobaro}
\frac{1}{2}\frac{d}{d} \intOB{ |\Div \, \vu|^2+ |\Grad\vr|^2}
     +  (2\mu+\lambda) \intO{ |\Grad \Div \, \vu|^2} = 0.
\end{equation}
 Thus \eqref{energybaro}, \eqref{entropybarofinal}, \eqref{curlbaro}, \eqref{divgradrhobaro}
provide an exponential decay in time for $\sqrt\kappa \|(\vr,\vu)\|_{H^1(\Omega)}$ to $0$ assuming initial perturbation $\sqrt\kappa (\vr_0,\vu_0)$ in $H^1(\Omega)$. Note in particular, that coefficient  $\kappa$ can be arbitrary small.

\subsection{Linearized heat-conducting compressible Navier--Stokes equations}\label{Linearized} In this subsection, we consider the linearized compressible Navier--Stokes equation with heat conductivity around $(1,0,1)$ namely
\begin{equation}\label{CNS}
\left\{
\begin{array}{l}
\partial_t\vr  + \Div \, \vu =0,\\
\displaystyle
\partial_t  \vu + \Grad \vr + \Grad \theta- 2\mu \Div  (D(\vu)) +\frac{2}{d} \mu \Grad \Div \, \vu  = \vc{0},\\
\displaystyle \partial_t \theta + \frac{2}{d} \Div \, \vu - K \Delta \theta = 0,
\end{array}
\right.
\end{equation}
with the periodic  boundary conditions.  Such system may be found for instance in \cite{Vi} page 51.   For this system, we easily check that there exists an energy, defined as 
\begin{equation}\label{energyheat}
\frac{1}{2} \Dt \intOB{ |\vu|^2 + |\vr|^2 + \frac{d}{2} |\theta|^2}
   + 2 \mu \intO{ \left|D(\vu)- \frac{\Div \, \vu}{d} \vc{I}\right|^2} + \frac{Kd}{2} \intO{ |\Grad\theta|^2} = 0.
\end{equation}
Differentiating the mass equation and the temperature equation with respect to 
space,  we find the following equations
\begin{equation}\label{NS-Aug}
\left\{
\begin{array}{l}
\partial_t\vr  + \Div \, \vu =0,\\
\displaystyle
\partial_t  \vu + \Grad \vr + \Grad \theta - 2\mu \Div  (D(\vu)) +\frac{2}{d} \mu \Grad \Div \, \vu  = \vc{0},\\
\displaystyle
\partial_t \theta + \frac{2}{d} \Div \, \vu- \kappa \lap \theta = 0, \\
\partial_t\Grad\vr + \Grad\Div \, \vu =\vc{0}, \\
\displaystyle
\partial_t\Grad\theta + \frac{2}{d} \Grad \Div \, \vu - \kappa \Delta \Grad\theta = 0. 
\end{array}
\right.
\end{equation}
Thus we can write a new equation satisfied by artificial velocity $\vu+ 2\kappa \mu\Grad\vr$
$$\partial_t  (\vu+2\kappa \mu\Grad\vr)
   + \Grad \vr + \Grad \theta 
   -  \mu \Delta \vu - \mu (1-\frac{2}{d}- 2 \kappa) \Grad \Div \, \vu  = 0.$$
Testing it against $\vu+2 \kappa \mu\Grad\vr$, we obtain
\eq{\label{entropyheat}
&\frac{1}{2} \Dt \intOB{ |\vu+ 2\kappa \mu \Grad\vr|^2
  + \vr^2} +2 \kappa \mu \intO{ |\Grad\vr|^2} + \mu \intO{ |\Grad \vu|^2}
+ \mu (1- \frac{2}{d} -  2\kappa) \intO{ |\Div \, \vu|^2 } \\
   & - 4\mu^2\kappa \lr{1-\frac{1}{d} - \kappa} \intO{ \Grad \Div \, \vu \cdot \Grad\vr}
   + \intO{\Grad\theta\cdot \vu} +2s \kappa \mu \intO{ \Grad\theta\cdot\Grad \vr} = 0.}
Thus as in the barotropic case, we have 
\eq{\label{entropyheat1}
&\frac{1}{2} \Dt \intOB{ |\vu+ 2 \kappa \mu \Grad\vr|^2
  + \vr^2} + 2 \kappa \mu \intO{ |\Grad\vr|^2} + 2 \mu \intO{ |A(\vu)|^2}
   + 2\mu (1 - \frac{1}{d} - \kappa) \intO{ |\Div \, \vu|^2 } \\
   & - 4 \mu^2\kappa \lr{1-\frac{1}{d} - \kappa} \intO{ \Grad \Div \, \vu \cdot \Grad\vr}
   + \intO{ \Grad\theta\cdot \vu} +  2\kappa \mu \intO{ \Grad\theta\cdot\Grad \vr }= 0.}
Recall that for $\Grad\vr$, we have
\begin{equation} \label{gradrhoheat}
\frac{1}{2} \Dt \intO{ |\Grad\vr|^2} + \intO{ \Grad \Div \, \vu \cdot \Grad \vr }= 0.
\end{equation}
Multiplying this relation by $4\displaystyle \mu^2\kappa \lr{1-\frac{1}{d} - \kappa}$ and adding to the previous one, we get
\eq{\label{entropyheat1}
&\frac{1}{2} \Dt \intOB{ |\vu+ 2 \kappa \mu \Grad\vr|^2 +  \vr^2
 + 4\mu^2\kappa \lr{1-\frac{1}{d} - \kappa}|\Grad\vr|^2} 
  + 2\kappa \mu \intO{ |\Grad\vr|^2 \\
  &+ 2 \mu \int |A(\vu)|^2}
   + 2 \mu \lr{1 - \frac{1}{d} - \kappa} \intO{ |\Div \, \vu|^2 }
   + \intO{ \Grad\theta\cdot \vu }+ 2\kappa \mu \intO{ \Grad\theta\cdot\Grad \vr} = 0.}
Recall that
\begin{equation}\label{tempp}
\frac{1}{2}\Dt \intO{ \frac{d}{2} |\theta|^2 }+ \frac{Kd}{2} \intO{ |\Grad\theta|^2} + \intO{  \theta \, \Div \, \vu }= 0,
\end{equation}
which, when added to \eqref{entropyheat1}, yields
\eq{\label{entropyheat2}
&\frac{1}{2} \Dt \intOB{ |\vu+ 2 \kappa \mu \Grad\vr|^2 +  \vr^2+ 4 \mu^2\kappa \lr{1-\frac{1}{d} - \kappa}|\Grad\vr|^2+\frac{d}{2} |\theta|^2} \\
    &+ 2 \mu \intO{ |A(\vu)|^2 }+ 2 \mu (1 - \frac{1}{d} - \kappa) \intO{ |\Div \, \vu|^2} \\
     &+ \kappa \mu \intO{ |\Grad\vr|^2} + \displaystyle (\frac{dK}{2}-\kappa\mu)\intO{ |\Grad\theta|^2}
   +  \kappa\mu\intO{ |\Grad(\vr+\theta)|^2} = 0.}
Choosing $\kappa$ such that $0<\kappa\mu <dK/2$ and $0<\kappa < (d-1)/d$, 
the $\kappa$-entropy balance  \eqref{entropyheat2} gives the exponential decay of $(\Grad\vr,\vu,\theta)$ in the $L^2$ norm, note in particular an interesting interplay between conductivity and pressure.

\bigskip
 Let us now focus on exponential decay in $H^1$ norms. 
 Let us take the equation satisfied by $\Div \, \vu$, test it against $\Div \, \vu$ and
integrate, we get
\begin{equation} \label{divheat}
\frac{1}{2} \Dt \intO{ |\Div \, \vu|^2} 
     + \intO{ \Delta \vr  \,  \Div \, \vu} + \intO{ \Delta \theta \, \Div \, \vu}
     + (2\mu+\lambda)\intO{ |\Grad \Div \, \vu|^2} = 0.
\end{equation}
 Let us now take the equation satisfied by $\Grad\theta$ and test it against
 $\Grad \theta$, we get
 \begin{equation}\label{gradthetaheat}
 \frac{1}{2}\Dt \intO{ |\Grad\theta|^2} - \frac{2}{d}\intO{ \Div \, \vu  \,  \Delta\theta}
     + K \intO{ |\Delta\theta|^2} = 0.
 \end{equation}
 Recall that 
 \begin{equation} \label{gradrhoheat1}
\frac{1}{2} \Dt \intO{ |\Grad\vr|^2} - \intO{ \Div \, \vu  \,  \Delta \vr }= 0.
\end{equation}
Thus adding \eqref{divheat}, \eqref{gradthetaheat} to \eqref{gradrhoheat1} gives
\begin{equation} \label{Estim}
\frac{1}{2} \Dt \intOB{ |\Div \, \vu|^2 
     +  |\Grad\vr|^2 
     +   \frac{d}{2} |\Grad\theta|^2 }
     + \frac{dK}{2}  \intO{ |\Delta\theta|^2}
     + (2\mu+\lambda)\intO{ |\Grad \Div \, \vu|^2} = 0.
\end{equation}
Recalling that the $\Curl\vu$ satisfies \eqref{curlbaro},
combining it with \eqref{Estim} and the $\kappa$-entropy equality \eqref{entropyheat2}, we get the exponential decay in time of  $\|\sqrt\kappa (\vr,\vu,\theta)\|_{H^1(\Omega)}$ if the initial
perturbation $\sqrt\kappa (\vr_0,\vu_0,\theta_0)$ is uniformly bounded in $H^1(\Omega)$.

\bigskip

\section{Heat-conducting Navier--Stokes equations with $\kappa$--energy.}
 In this section, we present the equations of motions for the heat-conducting fluid written in terms of the two velocities $\vu$ and $\vu+2\Grad\varphi(\vr)$ with corresponding densities $(1-\kappa)\vr$ and $\kappa\vr$.
We do not aim at proving the existence result for such system but on showing that the two-velocity hydrodynamics  in the spirit of the work by {\sc S.M. Shugrin} \cite{Sh} is consistent with the study performed for the low Mach number system in the first part of this diptych in \cite{BrGiZa}. More precisely, we will show that the formal low-Mach number limit for the two-velocities system gives the augmented system used in  \cite{BrGiZa} to construct the approximate solution.
An important observation is that the  system presented below does not coincide with the usual heat-conducting compressible Navier-Stokes equations. Indeed, the two-velocities description of the dynamics of the fluid lead to different energy equation with a generalised temperature, called {\it the $\kappa$-temperature}. However, this is not {\it a priori} the usual temperature, unless the system reduces to the angle velocity one (i.e. the density $\kappa\vr$ is equal to 0). This property was also explained in the works \cite{Sh} and \cite{GaSh} where the authors discuss {\it the capillary-temperature}. 
 
   We assume that $\Omega$ is a periodic box in $\R^3$, i.e. $\Omega=\mathbb{T}^3$, and we consider the following two-velocity system
        \begin{equation}\label{1.1}
     \left\{  
      \begin{array}{l}
            \vspace{0.2cm}
        \pt\vr+\Div (\vr \vu) = 0\\
            \vspace{0.2cm}
        \ptb{\vr\vu}+\Div (\vr \vu \otimes \vu) -  \Div(2 \mu(\vr) D(\vu)) - \Grad(\lambda(\vr) \Div u) + 
          \Grad p(\vr,e_{\kappa}) =\vc{0}\\
            \vspace{0.2cm}
        \pt(\vr(\vu+2\Grad\varphi(\vr)) + \Div (\vr \vu \otimes (\vu+2\Grad\varphi(\vr))) 
         -  \Div(2 \mu(\vr) A(u))   + \Grad p(\vr,e_{\kappa}) =\vc{0}\\
                \vspace{0.2cm}
        \pt{E_{\kappa}}+\Div (E_{\kappa}  \vu) + \Div(p\, [(1-\kappa)\vu + \kappa(\vu+2\Grad\varphi(\vr)]) 
           +\Div{\,\vQ_{\kappa}} \\
      \hskip2cm -\Div \Bigl((1-\kappa) \vS_1 (\vu) + \kappa\vS_2 (\vu+2\Grad\varphi(\vr))\Bigr)
        =0,\\
            \end{array}\right.
    \end{equation}
where we denoted  $D(\vu)=\frac{1}{2}\left(\Grad \vu+\Grad^t \vu\right)$ 
and $A(\vu) = \frac{1}{2}\left(\Grad \vu-\Grad^t \vu\right)$.\\
The viscosity coefficients $\mu(\vr)$, $\lambda(\vr)$ satisfy the  {\sc Bresch-Desjardins} relation
	\eq{\lambda(\vr)=2\mu'(\vr)\vr-2\mu(\vr).\label{BD1}}
  The total $\kappa$--energy $E_{\kappa}$ is defined as follows
\eq{E_{\kappa} = \vr \Bigl(e_{\kappa} +   \frac{(1-\kappa)}{2}|\vu|^2 
    + \frac{\kappa}{2}|\vu+2\Grad\varphi(\vr)|^2\Bigr). \label{Ekappa}}
\begin{rmk}
 Note that $E_\kappa \vu$ is expressed as a sum of two energies
$$(1-\kappa) \rho \vu \left[e_\kappa +\frac{|\vu|^2}{2}\right] +
                     \kappa \rho \vu  \left[e_\kappa +\frac{|\vu
                     +2\Grad\varphi(\vr)|^2}{2}\right]$$
 similarly to energy from {\rm \cite{Sh}}.
 \end{rmk}
 \begin{rmk}
Integrating the total $\kappa$--energy equation with respect to space, we obtain
$$\Dt\intO{E_\kappa}=0.$$
Thus \eqref{Ekappa} and  the identity 
 \begin{equation} \label{mixture}
    \frac{(1-\kappa)}{2}|\vu|^2 
    + \frac{\kappa}{2}|\vu+2\Grad\varphi(\vr)|^2=   \frac{|\vu + 2\kappa \Grad \varphi(\vr)|^2}{2} 
   + (1-\kappa)\kappa \frac{|2 \Grad\varphi(\vr)|^2}{2} 
\end{equation}
yields the following conservation property
$$
  \frac{d}{dt} 
    \intO {\vr \Bigl(e_{\kappa} + \frac{1}{2}|\vu+\kappa\Grad\varphi(\vr)|^2 
    + \frac{(1-\kappa)\kappa}{2}|\Grad\varphi(\vr)|^2\Bigr)} = 0.
$$
This quantity may be treated as a generalization  of the $\kappa$--entropy, found for the barotropic case, to the heat-conducting case.
\end{rmk}
The viscous tensors $\vS_1$, $\vS_2$ are given by
$$ 
   \vS_1 
   =  
   2\mu(\vr)D(\vu)+  \lambda(\vr)  \Div \, \vu \, \vc{I}
$$
and
$$  
  \vS_2 = 2\mu(\vr) A(\vu).
$$
The heat flux $\vc{Q}_\kappa$ is given by standard Fourier's law, i.e.
$$\vQ_\kappa=-K\Grad\theta_\kappa,$$
where $K$ is the positive heat-conductivity coefficient and $\theta_\kappa$ denotes the generalized temperature (the $\kappa$--temperature).
 Let us consider an ideal polytropic gas, namely
$$p(\vr, e_{\kappa}) = r\vr \theta_{\kappa} + p_c(\vr), \qquad 
    e_{\kappa} = C_v \theta_{\kappa} + e_c(\vr),$$
    where $r$ and $C_v$ are two constant positive coefficients, see for instance  \cite{BrDe3}.
For convenience, we denote $\gamma= 1+r/C_v$. Moreover, the additional
pressure and internal energy, $p_c$ and $e_c$ respectively, are associated to the ''zero Kelvin isothermal'', which roughly speaking means that 
$$\lim_{\vr\to0^+} p_c(\vr) =-\infty.$$
Further, we require that $e_c$ is a $C^2(0,\infty)$ nonnegative function and that the following constraint is satisfied 
$$p_c(\vr) = \vr^2 \frac{d e_c}{d\vr}(\vr).$$

Below we present two different formulations of the internal energy equation which lead to useful bounds on $\kappa$-temperature similarly as in \cite{BrDe3} for the usual temperature. \\
The first formulation reads
\eqh{
   & C_v\Bigl( \partial_t(\vr \theta_{\kappa}) + \Div (\vr \vu \theta_{\kappa}) 
     + \Gamma \vr \theta_{\kappa} \Div \, \vw\Bigr) \\
    & =  2(1-{\kappa})\mu(\vr)|D(\vu)|^2 + 2\kappa \mu(\vr)|A(\vu)|^2
        + 2(1-{\kappa})(\mu'(\vr)\vr-\mu(\vr))|\Div \, \vu|^2 + \Div(K\Grad \theta_{\kappa}),
}   
with $\Gamma$  the Gruneisen parameter and where the mixing temperature $\theta_{\kappa}$ becomes the usual temperature if $\kappa= 0$.
Note that for $0\le \kappa\le 1$, the $\kappa$-temperature remains non-negative in view of the maximum principle.\\
The second formulation is based on the notion of generalized  $\kappa$-entropy $s_{\kappa}$. It is the usual entropy in which the standard temperature has been replaced by the $\kappa$-temperature, i.e.
$$s_\kappa=C_v\log\theta_\kappa-r\log \vr,$$
thus, when $\vr,\theta_\kappa$ is sufficiently regular we can derive the following equation
 \eq{\label{entropy}
&\partial_t(\vr s_{\kappa})+ \Div(\vr\vu s_{\kappa})-\Div\lr{K\Grad\log\theta_\kappa}\\
&\qquad\qquad=   2(1-\kappa)\frac{\mu(\vr)|D(\vu)|^2}{  \theta_{\kappa} } + 2 \kappa\frac{\mu(\vr)|A(\vu)|^2}{  \theta_{\kappa} } + 2(1-\kappa)\frac{(\mu'(\vr)\vr-\mu(\vr))|\Div \, \vu|^2 }{  \theta_{\kappa} } \\
         &\qquad\qquad\quad- 2 \kappa \Gamma \vr \Delta\varphi(\vr)  + K\frac{|\Grad  \theta_{\kappa}|^2}{  \theta_{\kappa}^2}.
}
Note that, recalling the relation
 $$|D(\vu)|^2
     = \left|D(\vu) - \frac{1}{3}\Div  \vu\,  \vc{I}\right|^2
     + \frac{1}{3} |\Div  \vu|^2,$$
the terms on right-hand side, when integrated over space, give nonnegative contribution using
the assumption on $3\lambda(\vr) + 2\mu(\vr)$ and $\lambda(\vr) = 2(\mu'(\vr)\vr - \mu(\vr))$.
  Indeed, it suffices to check that for the penultimate term we have 
$$
  \intO {\frac{-\Gamma \vr \theta_{\kappa}  \Delta\varphi(\vr)}{\theta_{\kappa}}} 
     = \intO {\Gamma \varphi'(\vr)|\Grad\vr|^2} 
  \ge 0.
  $$
 Using all this information, it could be possible to prove global existence of $\kappa$-entropy
solution for the heat--conducting compressible-Navier Stokes system under analogous assumptions as in \cite{BrDe3}, replacing the usual temperature by the $\kappa$-temperature.  
The existence of the approximate solution could be then proven by using the augmented system written in terms of 
  $\vw= \vu+2\kappa\Grad\varphi(\vr)$ and $\vv= 2\Grad\varphi(\vr)$ as it was done in \cite{Sh}, or in Section \ref{S:constr} addressing  barotropic flows
\eq{\label{approxcond}
&\pt\vr+ \Div(\vr [\vw]_\delta) - \kappa\Div\lr{[\mu'( \vr)]_\alpha \Grad\vr} = 0,\\
&{}\\
& \pt\lr{\vr \vw}+ \Div  ((\vr [\vw]_{\delta} - \kappa  [\mu'( \vr)]_\alpha \Grad\vr) \otimes \vw)  
 -  \Grad (\lr{\lambda(\vr)- \kappa (\mu'(\vr)\vr-\mu(\vr))}\Div(\vw-\kappa\vv)) \\
  &   - (2-\kappa)  \Div (\mu(\vr) D(\vw)) - \kappa \Div (\mu(\vr) A(\vw))  
  + \varepsilon \lap^{2s} \vw - \varepsilon \Div ((1+|\Grad\vw|^2)\Grad\vw) + \Grad p(\vr,e_\kappa)\\
  &= - \kappa (2-\kappa) \Div(\mu(\vr) \Grad\vv),\\
&{}\\
&\partial_t(\vr \vv) + \Div((\vr \vw- \kappa [\mu'( \vr)]_\alpha \Grad\vr)\otimes \vv) 
- \kappa \Div(\mu(\vr)\Grad\vv) +\Grad ((\mu'(\vr)\vr - \mu(\vr))\Div \lr{\vw-\kappa\vv})\\
&= - \Div(\mu(\vr) \Grad^{t} \vw)
} 
with the $\kappa$-total energy supplemented by the $\ep$ correction corresponding to the $\ep$ regularisation of the momentum
\eqh{
&\pt{(\vr E_\kappa)}+\Div (\bigl(\vr [\vw]_\delta - \kappa  [\mu'( \vr)]_\alpha\Grad\vr \bigr) E_\kappa)
          +\Div(p\, \vw) 
          +\Div{\vQ}   \\
          &-\Div \Bigl(\vS_1\vw  +  (2-\kappa)\kappa \, \vS_2 \vv\Bigr)
          + \varepsilon \Delta^s \vw \cdot \vw - \varepsilon \Div ((1+|\Grad\vw|^2)|\vw|^2 =0
}
and the set of initial conditions. 
Above, the  total $\kappa$--energy $E_\kappa$ is defined as
$$E_\kappa = e_\kappa + \frac{1}{2}|\vw|^2 + \frac{\kappa(2-\kappa)}{2}|\vv|^2.$$
Note, however,  this construction would not lead to the usual heat--conducting compressible Navier-Stokes system in the limit $\ep\to0$. Indeed, the difference is again due to $\kappa$-temperature that is not the usual one. But, performing  a formal low Mach number limit for this system, we would get $p=1$, $\Div \, \vw=0$ 
(comparing terms of the same order).  In the equation on $\vw$, being now incompressible, the pressure gradient $\nabla p$ would be replaced by Lagrangian multiplier $\nabla\pi$. As a result, we would get the augmented system defined in \cite{BrGiZa} in the part devoted to construction of solution.
\bigskip

\noindent {\bf Acknowledgments.}
The authors want to  thank the referee for his/her valuable comments which allows to improve the quality of the paper.
   The first author acknowledges support from the ANR-13-BS01-0003-01 project DYFICOLTI.
The authors want to thank {\sc S.L. Gavrilyuk} for really interesting comments and for pointing to them the two interesting papers linked to media  with equations of state that depend on derivatives. The third author acknowledges the Post-Doctoral support of Ecole Polytechnique, she was also supported by MN grant  IdPlus2011/000661 and by the fellowship START of the Foundation for Polish Science.

%%%%%%%%%%%%%%%

\end{document}